\newif\iflatexe
\latexefalse

\iflatexe
\documentclass[11pt]{article}
\usepackage{latexsym,mathptm}
\else
\documentstyle[11pt]{article}
\fi
\makeatletter
\@addtoreset{equation}{subsection}

\makeatother

\newcommand{\separate}{\medskip\noindent}

\def\newsection{ \separate
   \refstepcounter{subsection} 
   {\large\bf \thesubsection\kern.3em}
}
\def\mytheorem#1{
   \separate{\large\bf Theorem#1:\kern.3em}} 

\def\mylemma#1{
   \separate{\large\bf Lemma#1:\kern.3em}} 

\def\mycorollary#1{
   \separate{\large\bf Corollary#1:\kern.3em}} 

\def\myprop#1{
   \separate{\large\bf Proposition#1:\kern.3em}} 

\newcommand{\GANZ}{{\sf Z\hspace*{-0.4em}Z}}
\newcommand{\REELL}{{\setlength{\unitlength}{1em}
                     \begin{picture}(0.75,1)
                     \put(0,0){\line(0,1){0.69}}
                     \put(0,0){\sf R}
                     \end{picture}
                   }}
\newcommand{\KOMPLEX}{{\setlength{\unitlength}{1em}
                     \begin{picture}(0.7,1)
                     \put(0.34,0){\line(0,1){0.65}}
                     \put(0,0){\sf C}
                     \end{picture}
                   }}
\newcommand{\NATUR}{{\setlength{\unitlength}{1em}
                     \begin{picture}(0.75,1)
                     \put(0,0){\line(0,1){0.69}}
                     \put(0,0){\sf N}
                     \end{picture}
                   }}
\newcommand{\FGANZ}{\mbox{\tiny{\rm Z\hspace*{-0.45em}Z}}}
\newcommand{\FNATUR}{\mbox{\tiny{\setlength{\unitlength}{1em}
                     \begin{picture}(0.6,0.5)
                     \put(0,0){\line(0,1){0.48}}
                     \put(0,0){\rm N}
                     \end{picture}
                   }}}
\newcommand{\FKOMPLEX}{\mbox{\tiny{\setlength{\unitlength}{1em}
                               \begin{picture}(0.6,0.5)
                               \put(0.34,0){\line(0,1){0.47}}
                               \put(0,0){\rm C}
                               \end{picture}
                              }}}
\newcommand{\FREELL}{\mbox{\tiny{\setlength{\unitlength}{1em}
                     \begin{picture}(0.6,0.5)
                     \put(0,0){\line(0,1){0.48}}
                     \put(0,0){\rm R}
                     \end{picture}
                   }}}
\newcommand{\Brr}{{\mathchoice{\REELL}{\REELL}{\!\FREELL}{\!\FREELL}}}
\newcommand{\Bcc}{{\mathchoice{\KOMPLEX}{\KOMPLEX}{\!\FKOMPLEX}{
\!\FKOMPLEX}}}
\newcommand{\Bnn}{{\mathchoice{\NATUR}{\NATUR}{\FNATUR}{\FNATUR}}}
\newcommand{\Bii}{{\mathchoice{\GANZ}{\GANZ}{\FGANZ}{\FGANZ}}}

\def\QED{\hfill$\Box$}
\def\inv{^{-1}}
\def\BEA{\begin{eqnarray}}
\def\EEA{\end{eqnarray}}
\def\BEQ{\begin{equation}}
\def\EEQ{\end{equation}}
\def\bref#1{(\ref{#1})}
\def\tmatrix#1#2#3#4{
    \left(\begin{array}{cc} #1 & #2 \\ #3 & #4 \end{array}\right)}
\def\LieSL{{\mbox{\bf SL}}}
\def\LieSO{{\mbox{\bf SO}}}
\def\LieGL{{\mbox{\bf GL}}}
\def\LieSU{{\mbox{\bf SU}}}
\def\Liesl{{\mbox{\bf sl}}}
\def\Liesu{{\mbox{\bf su}}}
\def\dprime{{\prime\prime}}
\def\TD{{T_D}}
\def\TU{{T_U}}
\def\TV{{T_V}}
\def\vecr{{\mbox{\bf r}}}
\def\vecsigma{{\sigma\kern-2.25mm\sigma}}
\def\ta{\tilde{a}}
\def\tb{\tilde{b}}
\def\tc{\tilde{c}}

\def\tf{\tilde{f}}
\def\tg{\tilde{g}}
\def\txi{\tilde{\xi}}
\def\tF{\tilde{F}}
\def\hf{\hat{f}}
\def\hg{\hat{g}}
\def\hxi{\hat{\xi}}
\def\hA{\hat{A}}
\def\FA{{\cal A}}
\def\FM{{\cal M}}
\def\FU{{\cal U}}
\def\der#1#2{{{\rm d}#1\over {\rm d}#2}}

\def\diff{{\rm d}}
\def\dz{\diff z}
\def\dF{\diff F}
\def\dhg{\diff\hg}
\def\dtg{\diff\tg}
\def\dzbar{\diff\zbar}
\def\dg{\diff g}
\def\Ad{{\mbox{\rm Ad}}}
\def\res{{\mbox{\rm res}}}
\def\diag{{\mbox{\rm diag}}}

\def\Aut{{\mbox{\rm Aut}}}
\def\tr{{\mbox{\rm tr}}}
\def\const{{\mbox{\rm const}}}

\def\hg{\hat{g}}

\def\hF{\hat{F}}
\def\hp{\hat{p}}
\def\ts{\tilde{s}}
\def\hU{\hat{U}}
\def\hV{\hat{V}}

\def\tg{\tilde{g}}
\def\ty{\tilde{y}}

\def\tQ{\tilde{Q}}
\def\zbar{{\bar{z}}}

\def\Qbar{{\bar{Q}}}

\begin{document}

\renewcommand{\thefootnote}{\fnsymbol{footnote}}
\begin{center}
{\LARGE Meromorphic Potentials and Smooth CMC Surfaces}

\vskip1cm
\begin{minipage}{6cm}
\begin{center}
J.~Dorfmeister\footnotemark[1]\footnotemark[2]\\ Mathematisches Institut\\TU
M\"unchen\\D-80333 M\"unchen
\end{center}
\end{minipage}
\begin{minipage}{6cm}
\begin{center}
G.~Haak\footnotemark[3]\\ Department of Mathematics\\ University of Kansas\\
Lawrence, KS 66045
\end{center}
\end{minipage}
\vspace{0.5cm}

\end{center}

\footnotetext[1]{partially supported by NSF Grant 
DMS-9205293 and Deutsche For\-schungs\-ge\-mein\-schaft}
\footnotetext[2]{permanent address: Department of Mathematics,
University of Kansas, Lawrence, KS 66045}
\footnotetext[3]{supported by KITCS grant OSR-9255223}

\section{Introduction}\refstepcounter{subsection}
In the recent paper \cite{DPW} it was shown that every conformal CMC immersion
$\Phi:D\rightarrow \Brr^3$, $D\subset\Bcc$ the whole complex plane or the open
unit disk, can be produced from a meromorphic matrix valued one form
\BEQ \label{MEROPOT}
\xi=\lambda\inv\tmatrix{0}{f(z)}{g(z)}{0}\dz,
\EEQ
$\lambda\in S^1$, the so called ``meromorphic potential''.

Here $f$ and $g$ are meromorphic functions of $z\in D$ and
$f(z)g(z)=E(z)$, where $E(z)\dz^2$ is up to a constant factor (see the
appendix) the Hopf differential of the 
surface $M=(D,\Phi)$.

Moreover, the normalizations used in \cite{DPW} imply, that $f$ and
$g$ don't have any poles at $z=0$.

\noindent The construction involves the following steps:
\begin{enumerate}
\item Solve the initial value problem
\BEQ
\diff g_-=g_-\xi,\, g_-(z,\lambda)\in\Lambda^-_\ast\LieSL(2,\Bcc)_\sigma,\; 
g_-(0,\lambda)=I.
\label{ONE}
\EEQ
\item Compute an Iwasawa decomposition of $g_-$:
\BEQ
g_-=F g_+\inv,\,F\in\Lambda\LieSU(2)_\sigma,\;
g_+\in\Lambda^+\LieSL(2,\Bcc)_\sigma, \label{TWO}
\EEQ
where $g_+(\lambda=0)$ is of the form $\diag(a,a\inv)$ with $a$ an
arbitrary complex number as
opposed to the traditional Iwasawa decomposition, where $|a|=1$.
\item Apply the Sym-Bobenko formula,
\BEQ \label{SymBob}
\phi_\lambda(z)=\der{}{\theta}F\cdot F\inv
+{1\over2}F\tmatrix{i}{0}{0}{-i}F\inv,\;
\lambda=e^{i\theta},
\EEQ
to $F=F(e^{i\theta},z)$ to obtain CMC immersions with mean 
curvature~$H=-{1\over2}$.
This actually yields a CMC immersion for every $\lambda\in S^1$ 
(see~\ref{Sym}). The
immersion $\Phi=\Phi_1$ is the given one.
\end{enumerate}
Here $\Lambda\LieSL(2,\Bcc)_\sigma$ and $\Lambda\LieSU(2)_\sigma$
denote the twisted loop groups over $\LieSL(2,\Bcc)$ and
$\LieSU(2)$, respectively, given by the automorphism
$$
\sigma: g\mapsto(\mbox{\rm Ad}\sigma_3)(g),\,\sigma_3=\tmatrix{1}{0}{0}{-1},
$$
e.g.\ 
$$
\Lambda\LieSL(2,\Bcc)_\sigma=\{g:S^1\longrightarrow\LieSL(2,\Bcc)|\;
g(-\lambda)=\sigma(g(\lambda))\}.
$$
In order to make these loop groups complex Banach Lie groups, we equip them,
as in \cite{DPW}, with some $H^s$-topology for $s>{1\over2}$.

Elements of these twisted loop
groups are matrices with offdiagonal entries which are odd functions and
diagonal entries which are even functions in the parameter $\lambda$.
All entries are in the Banach algebra ${\cal A}$ of $H^s$-smooth functions.

Their Lie-algebras are then complex Banach Lie algebras, e.g.\
$$
\Lambda\Liesl(2,\Bcc)_\sigma=\{x: S^1\longrightarrow\Liesl(2,\Bcc)|\;
x(-\lambda)=\sigma(x(\lambda)),\; \mbox{\rm $x$ is $H^s$-smooth}\}.
$$
The indices $+$ and $-$ refer to the usual splitting of the loop group
into maps analytic inside and outside the unit circle.
In addition by $\Lambda^-_\ast\LieSL(2,\Bcc)$ we denote the subgroup
of elements of $\Lambda^-\LieSL(2,\Bcc)$, that take the
value $I$ at infinity. Then the set
$\Lambda^-_\ast\LieSL(2,\Bcc)_\sigma\cdot\Lambda^+\LieSL(2,\Bcc)_\sigma$ 
is open and
dense in $\Lambda\LieSL(2,\Bcc)_\sigma$. 

The decomposition, 
$$
\Lambda\LieSL(2,\Bcc)_\sigma\cong\Lambda\LieSU(2)_\sigma\times
\Lambda^+\LieSL(2,\Bcc)_\sigma,
$$
is defined for all elements of $\Lambda\LieSL(2,\Bcc)_\sigma$. 

To be more precise, in \cite{DPW} the following twisted versions of
well known theorems \cite{PS} are proved:

\separate (i) For each solvable subgroup $B$ of $\LieSL(2,\Bcc)$ with
$\LieSU(2)\cap B=\{I\}$, multiplication 
$$
\Lambda\LieSU(2)_\sigma\times\Lambda^+_B\LieSL(2,\Bcc)_\sigma
\longrightarrow\Lambda\LieSL(2,\Bcc)_\sigma
$$
is a diffeomorphism onto. Here $\Lambda^+_B\LieSL(2,\Bcc)_\sigma$ denotes
the set of all elements $g_+(\lambda)$ of
$\Lambda^+\LieSL(2,\Bcc)_\sigma$ with $g_+(0)\in B$.

\separate (ii) Multiplication 
$$
\Lambda^-_\ast\LieSL(2,\Bcc)_\sigma\times\Lambda^+\LieSL(2,\Bcc)_\sigma
\longrightarrow\Lambda\LieSL(2,\Bcc)_\sigma
$$
is a diffeomorphism onto the open and dense subset 
$\Lambda^-_\ast\LieSL(2,\Bcc)_\sigma\cdot\Lambda\LieSL^+(2,\Bcc)_\sigma$
of $\Lambda\LieSL(2,\Bcc)_\sigma$, called the ``big cell'' \cite{SW86}.

If at a point $z_0\in D$, $f$ has a zero and $g$ is holomorphic, then,
as was pointed out in \cite{DPW} the map $\Phi_1$ fails to be an
immersion at $z=z_0$ (see also the appendix). 

If $g(z)$ has a zero at $z=z_0$, and $f$ is defined and nonzero there,
then one gets an umbilic at the image of $z_0$ in
$\Brr^3$.

This allows one to use the DPW method, to construct a CMC immersion
with a prescribed distribution of umbilics in the domain $D$. 

Special cases involve the sphere (minus a point), 
where $E(z)=0$, $f=1$,
the cylinder with $E(z)=f(z)= 1$ and the generalized Smyth
surfaces, where $E(z)$ is a polynomial whose roots give the umbilics,
and $f(z)= 1$. 

These cases are among many investigated and visualized by several
groups using the Bubbleman AVS network by Ulrich Pinkall and Charlie
Gunn, which is based on earlier work of D. Lerner and I. Sterling~\cite{LS}. 
In all cases the functions $f(z)$ and $g(z)$ were
chosen holomorphic in the entire domain $D$.

If we start with a meromorphic function $f(z)$ and a nonvanishing 
holomorphic Hopf
differential $E(z)\dz^2$, each of the three steps in the construction
imposes conditions on the choice of the meromorphic potential, i.e.\
on $f(z)$ and $E(z)$.  In order to get a smooth surface, we first need
to find a meromorphic solution
$g_-(z,\lambda)\in\Lambda_\ast^-\LieSL(2,\Bcc)_\sigma$ to \bref{ONE},
which will turn out not to be possible for arbitrary functions
$f(z)$ and $E(z)$.

In the second step the factor $F$ in the Iwasawa 
decomposition needs to be
smooth, which imposes further conditions on $f(z)$ and $E(z)$.

Finally, applying the Sym-Bobenko formula (see the appendix), we
obtain a map $\Phi=\Phi_1$, but one can get branch points of the
surface defined by $\Phi$, i.e.\ $\Phi=\Phi_1$ possibly fails to be an
immersion at certain points of $D$.

In this note we will give necessary and sufficient conditions on
$f(z)$ and $E(z)$ for $\Phi$, the CMC map associated with the
meromorphic potential $\xi$, to be a CMC immersion. Then the
surface $M=(D,\Phi)$ will be smooth without branchpoints. On the other
hand, if $M$ is a smooth immersed CMC surface without branchpoints, then
there always exists a CMC immersion $\Phi$ and a subset $D$ of $\Bcc$, 
such that $M=(D,\Phi)$. Therefore, we locally describe all smooth
immersed CMC surfaces in $\Brr^3$ without branchpoints. 

In section~\ref{firstpart} we will develop a necessary and sufficient
algebraic condition for $f$ and the Hopf differential, which ensures
that there exists a meromorphic matrix solution of \bref{ONE}
in the proper twisted loop group.

While this gives us a meromorphic mapping $g_-(z,\lambda)$ from $D$
into the loop group $\Lambda^-_\ast\LieSL(2,\Bcc)_\sigma$, 
it does not guarantee
the smoothness of the compact part $F(z,\lambda)$ of the Iwasawa
decomposition \bref{TWO} of $g_-(z,\lambda)$.

Necessary and sufficient local conditions for this are derived in
section~\ref{secondpart}. These take the shape of compatibility
conditions on the pole and zero orders of the functions $f(z)$ and $E(z)$.
Incorporated are also conditions for the (non)existence of branch points.

In all our investigations we will exclude the trivial case $E=0$, the
case of the round sphere.

In section~\ref{examples} we will give three examples of
nonholomorphic meromorphic potentials, which are associated with
CMC immersions. 

We close with an appendix, in which we fix
the conventions used in this paper.

The authors would like to thank David Lerner for his continuing
interest and helpful information on the Sym-Bobenko formula, as well as Fran
Burstall and Franz Pedit for interesting discussions on the subject.

Most of this paper was written during a visit of one of the authors
(J.~D.) at the TU-M\"unchen. He would like to thank the TU-M\"unchen
for its hospitality.

The pictures at the end of the paper were produced using a
modification of the Bubbleman AVS network written by Charlie Gunn at
the SFB-288, TU-Berlin. 

\section{Meromorphic Solutions of $g_-\inv\diff g_-=\xi$} \label{firstpart}

\newsection \label{sec21}
In this section we want to find a 
meromorphic matrix solution
$g_-(z,\lambda)$ of the initial value problem
\BEA
g_-^\prime=g_-\xi,\label{THREEA}\\
g_-(0,\lambda)=I.\label{THREEB}
\EEA
Here $(\cdot)^\prime$ denotes the derivative w.r.t.\ the complex
variable $z$ and
\BEQ
\xi=\lambda\inv Q(z)=\lambda\inv\tmatrix{0}{f}{g=f\inv E}{0},
\EEQ
i.e.\ we have omitted the $\dz$ factor in equation \bref{ONE}.

We note that then automatically 
$g_-(z,\lambda)\in\Lambda^-_\ast\LieSL(2,\Bcc)_\sigma$, i.e.\
$\det g_-(z,\lambda)= 1$, the diagonal entries are even
functions and the offdiagonal entries are odd
functions of $\lambda$, and $g_-(z,\infty)= I$.

It is important to note that the equation \bref{THREEA}
is invariant under coordinate changes. More precisely, if $f$ and $E$
are the defining data for $Q$ in the coordinate $z$, then in
a different coordinate $w$ the corresponding functions $\hat{f}(w)$ 
and $\hat{E}(w)$ are given by
\BEA \label{fchange}
\hat{f}(w) & = & f(z(w))\der{z}{w},\\
\hat{E}(w) & = & E(z(w))\left(\der{z}{w}\right)^2.
\EEA
If in addition $w=0\Leftrightarrow z=0$, then also the equation
\bref{THREEB} is valid in both coordinate systems.

This allows the following normalization of $f$: 
Locally around a point $z_0$ we can always choose the
coordinate $w$ to be such that $\hat{f}(w)=(w-w_0)^k$ or 
$\hat{f}(w)=(w-w_0)^{-k}$, where
$w_0=w(z_0)$ and $k\geq0$ is the pole/zero-order of $f$ at $z_0$.

Next we look at the matrix entries of $g_-$:
\BEQ
g_-=\tmatrix{a}{b}{c}{d}.
\EEQ
Equation \bref{THREEA} can be written as a set of four scalar
differential equations:
\BEA
\lambda a^\prime & = & b f\inv E,\label{FOURA}\\
\lambda b^\prime & = & a f,\label{FOURB}\\
\lambda c^\prime & = & d f\inv E,\label{FOURC}\\
\lambda d^\prime & = & c f,\label{FOURD}
\EEA
and the initial condition \bref{THREEB} translates to
\BEA
a(0,\lambda)=d(0,\lambda)& = & 1\label{FIFTEEN},\\
b(0,\lambda)=c(0,\lambda) & = & 0.\label{SIXTEEN}
\EEA
We differentiate equations \bref{FOURB} and \bref{FOURD} and
substitute for $a^\prime$ and $c^\prime$ using equations \bref{FOURA}
and \bref{FOURC}, respectively. 
It follows that $b$ and $d$ solve 
\BEQ
y^\dprime-{f^\prime\over f}y^\prime-\lambda^{-2}E y=0. \label{FIVE}
\EEQ
Analogously we show that $a$ and $c$ solve 
\BEQ
y^\dprime+\left({f^\prime\over f}-
{E^\prime\over E}\right)y^\prime-\lambda^{-2}E y=0. \label{SIX}
\EEQ
The initial conditions are given by \bref{FIFTEEN} and \bref{SIXTEEN}.

\newsection
Note that \bref{SIX} can also be stated as 
\BEQ \label{SIXA}
y^\dprime-{g^\prime\over g}y^\prime-\lambda^{-2}E y=0.
\EEQ

\separate{\large\bf Remark:} In this paper we will use methods and
results of \cite{DPW}. Therefore we need that all functions of
$\lambda$ are in the algebra $\FA=H^s$, $s>{1\over2}$. We note that
the coefficients of the
differential equations above are, as functions of $\lambda$, obviously
contained in the algebra $\FA$. Therefore by the standard theory of
differential equations in Banach spaces (see e.g.~\cite{Deim}) the
solutions $y=y(z,\lambda)$ are, as functions of $\lambda$, analytic
and contained in $\FA$ as well. Moreover, if we assume
$y(0,\lambda)=1$, then $y$ has an expansion relative to $\lambda$ of
the form
\BEQ
y=\sum_{n=0}^\infty q_n\lambda^{-2n}.
\EEQ 
Here the coefficients are functions of $z$. If $y$ is meromorphic in
$z$, then the $q_n=q_n(z)$ are meromorphic as well. So all the
functions occurring in this paper will have values in the Banach
algebra $\FA=H^s$, $s>{1\over2}$, as required by \cite{DPW}.

\label{th21}
\mytheorem{} {\em Let $f$ and $g$ be meromorphic
functions without poles at $z=0$.
The initial value problem \bref{THREEA}, \bref{THREEB} with
meromorphic potential 
\bref{MEROPOT} has a global meromorphic matrix solution $g_-$ if and
only if \bref{FIVE} has at every point in $D$ two linearly independent
local meromorphic solutions.}

\separate
We will postpone the proof of this theorem until section~\ref{proof21}.

\separate{\large\bf Remark:} 
In view of equations \bref{FOURB} and \bref{FOURD} we note that every
solution of \bref{FIVE} induces the solution $f\inv y^\prime$ of 
\bref{SIX}. Similarly every solution of \bref{SIX} produces a
solution of \bref{FIVE}. Moreover the
existence of two linearly independent meromorphic solutions of
\bref{FIVE} is equivalent to the existence of two meromorphic
solutions of \bref{SIX}.

\newsection Theorem~\ref{th21} shows, that it suffices to 
investigate \bref{FIVE}.

If $f$ has neither poles nor zeroes around $z=z_0$, the second 
order differential equation \bref{FIVE} has locally holomorphic 
coefficients. Therefore there exist locally two linearly independent 
holomorphic solutions.
So let us look at the poles and zeroes of $f$.

We will first prove a theorem which allows us to restrict our
attention to poles and zeroes of even order.

\mytheorem{} {\em The function $f$ in the meromorphic
potential is the square of a meromorphic function.}

\separate\underline{Proof:}
By \cite{DMPW} there exists for 
a given CMC immersion $\Phi$ a holomorphic potential $\eta$ of the
form
\BEQ
\eta(z,\lambda)=\lambda\inv\tmatrix{0}{c}{c\inv
E(z)}{0}+\eta_+(z,\lambda),\, c\in\Bcc,\;c\neq0
\EEQ
where $\eta_+(z,\lambda)\in\Lambda^+\Liesl(n,\Bcc)_\sigma$.
The Birkhoff splitting $g=g_-g_+$ of the integral $g$ of $\eta$
produces the meromorphic potential. But this amounts to a gauge
transformation of $\eta$ with $g_+\inv$. As $g_+$ contains no negative
powers of $\lambda$, the $\lambda\inv$ coefficient of
$g_-\inv\dg_-$ is the $\lambda\inv$ coefficient of $\eta$
conjugated with the constant term $g_0$ of $g$. If we write
$g_0=\diag(\omega_0,\omega_0\inv)$, then the meromorphic potential is given by
\BEQ
\xi=g_-\inv\dg_-=g_+\eta
g_+\inv-\dg_+g_+\inv=\lambda\inv\tmatrix{0}{cw_0^2}{c\inv w_0^{-2}E}{0}\dz.
\EEQ
and $f=cw_0^2$. \QED

\separate It follows directly:

\mycorollary{} {\em $f$ has only poles and zeroes of even order.}

\newsection \label{sec23}
Assume $f$ has a pole of order $n\geq2$ at some point $z=z_0$.
Then \bref{FIVE} has a regular singular point at $z=z_0$ and the
form of the solutions can be decided by looking at the
indicial equation 
\BEQ
r(r-1)+nr=0 \label{SEVEN}
\EEQ
The roots of this equation are
\BEQ
r_1=0,\,r_2=-n+1. \label{EIGHT}
\EEQ
The theory of ordinary differential equations with regular singular
points (see e.g.\
\cite{C}) shows that there always exists a locally meromorphic
solution $y_1$ of \bref{FIVE}. This solution is associated with the
higher of the two roots.
It is also well known, that every formal 
power series solution of \bref{FIVE} converges.
 
Making the ansatz $y_2=y_1\cdot v$ for a second linearly independent 
solution $y_2$, we get
\BEQ
v^\dprime = \left(\ln{f\over y_1^2}\right)^\prime v^\prime,
\EEQ
or
\BEQ \label{resorig}
v^\prime=C{f\over y_1^2},\;C=\const\neq0
\EEQ
In general this second solution will therefore have a
logarithmic singularity at $z=z_0$, unless the residue of the right
hand side of \bref{resorig} vanishes.
\BEQ
\oint_{z_0}{f\over y_1^2}\diff z=\res_{z=z_0}{f\over y_1^2}=0. \label{TEN}
\EEQ
The solution $y_1$ can be written as a power series ($w=z-z_0$)
\BEQ \label{y1caseone}
y_1=\sum_{i\geq0} a_iw^i,\, a_0=1, \label{NINE}
\EEQ
which has has neither a pole nor a zero at $z_0$. $y_1$ is a
nontrivial holomorphic function near $z=z_0$.

\newsection \label{sec24}
Let us now assume that $f$ has a zero of order $n\geq2$ at $z=z_0$.
This case is only interesting
if $f$ doesn't divide the Hopf differential $E$, because otherwise $g$ is also
nonsingular at $z=z_0$, and we obtain two holomorphic solutions of \bref{FIVE}.

Then in analogy to the last section, one gets the indicial equation
\BEQ
r(r-1)-nr=0,
\EEQ
with roots 
\BEQ
r_1=n+1,\,r_2=0. \label{EIGHT2}
\EEQ
We denote by $y_1$ the solution associated with the larger root
$r_1$. It is locally holomorphic and of the form ($w=z-z_0$)
\BEQ \label{y1casetwo}
y_1=\sum_{i\geq0}a_iw^{i+n+1},\;a_0=1.
\EEQ
The ansatz $y_2=y_1\cdot v$ again yields condition \bref{TEN} for the
second solution $y_2$ to be meromorphic. $y_2$ is then also locally
holomorphic at $z=z_0$.

\mytheorem{} {\em If $f$ has a pole or zero at $z=z_0$ and if $y_1$ denotes
a meromorphic solution to \bref{FIVE}, then the following is
equivalent:
\begin{description}
\item[a)] \bref{FIVE} has locally two linearly independent meromorphic
solutions.
\item[b)] The solution $y_1$ satisfies equation \bref{TEN}.
\end{description}}

\separate\underline{Proof:} By the calculations in the last two sections
we know, that for every meromorphic solutions $y_1$, there is a second
solution $y_2$, given by $y_2=y_1\cdot v$, where 
$v$ satisfies \bref{resorig}. This solution is linearly independent of
$y_1$, since $f$ and therefore $v$ and $v^\prime$ don't vanish identically.
Moreover, $y_2$ is also meromorphic, if and
only if \bref{TEN} is satisfied. This proves the theorem.
\QED

\separate{\large\bf Remark:} By equation \bref{fchange} condition
\bref{TEN} is invariant under coordinate changes $z\rightarrow w(z)$.

\newsection \label{recursionsec}
The condition \bref{TEN} for the existence of two linearly independent
meromorphic solutions of \bref{FIVE} can be expressed more explicitly
in terms of $E$ and $f$.

In general $y_2$ will be of the form ($w=z-z_0$)
\BEQ
y_2(w)=\alpha y_1(w)\ln w+w^{r_2}\left(1+\sum_{k=1}^\infty \ta_kw^k\right).
\EEQ
If \bref{TEN} is satisfied, then the coefficient $\alpha$ of the
logarithmic term vanishes and $y_2$
is of the form
\BEQ
y_2=w^{r_2}+\sum_{k=1}^\infty \ta_k w^{k+r_2}.\label{TWFIVE}
\EEQ
We would also like to note that, if \bref{FIVE} has a solution of the
form \bref{TWFIVE} even only formally, where $r_2$ is the smaller of
the two roots of the indicial equation, then locally there exist two 
linearly independent meromorphic solutions.

Substituting \bref{TWFIVE} into \bref{FIVE} yields the recursion relation
($k\geq2$)
\BEA \label{recrel}
\varphi(k+r_2)\ta_k & = & -\sum_{s=0}^{k-1}(s+r_2)q_{k-s}\ta_s+\lambda^{-2}
\sum_{s=2}^k E_{s-2}\ta_{k-s}, \nonumber\\
\varphi(r_2+1)\ta_1 & = & -r_2 q_1 \ta_0,
\EEA
for the coefficients $\ta_k$, where $\varphi$ is the left hand side of
the indicial equation, i.e.\
\BEQ
\varphi(m)=m(m-1)+nm,
\EEQ
if $z=z_0$ is a pole of $f$ of order $n$, and
\BEQ
\varphi(m)=m(m-1)-nm,
\EEQ
if $z=z_0$ is a zero of $f$ of order $n$. 

Moreover, the coefficients $q_j$ are defined by 
\BEQ
-w\frac{f^\prime(w)}{f(w)}=\sum_{j=0}^\infty q_j w^j
\EEQ
and 
\BEQ
E(w)=\sum_{k=0}^\infty E_k w^k.
\EEQ
It is clear, that \bref{recrel} defines $\ta_k$ uniquely from
$\ta_0,\ldots,\ta_{k-1}$, as long as $\varphi(k+r_2)\neq0$. If we
choose $\ta_0=1$, we get by \bref{recrel}, 
that all $\ta_k$ are polynomials in $\lambda^{-2}$.
Therefore:

\mytheorem{} {\em
Let $r_1>r_2$ be the two solutions to the indicial equation of
\bref{FIVE}. The following statements are equivalent
\begin{enumerate}
\item The equation \bref{FIVE} has a solution of the form
\bref{TWFIVE}.
\item \bref{FIVE} has two linearly independent meromorphic solutions.
\item 
\BEQ \label{recrel2}
\sum_{s=0}^{r_1-r_2-1}
(s+r_2)q_{r_1-r_2-s}\ta_s=
\lambda^{-2}\sum_{s=2}^{r_1-r_2}E_{s-2}\ta_{r_1-r_2-s},
\EEQ
where $-w\frac{f^\prime(w)}{f(w)}=\sum_{j=0}^\infty q_j w^j$ and 
$E(w)=\sum_{k=0}^\infty E_k w^k$.
\end{enumerate}
Moreover, $r_1=0$, $r_2=-n+1$, $q_0=n$, if $z=z_0$ is a pole of $f$ of order
$n$, and $r_1=n+1$, $r_2=0$, $q_0=-n$, 
if $z=z_0$ is a zero of $f$ of order $n$.}

\separate{\large\bf Remark:} 
Condition~\bref{recrel2} above is independent of the special shape of
the function $f$. In sections~\ref{normalized}--\ref{lastnorm} 
we will further investigate it for
the local normalizations $f(z)=(z-z_0)^n$ and $f(z)=(z-z_0)^{-n}$.
In these cases it is possible to derive a condition on the functions
$f$ and $E$, which is more explicit and computationally easier to handle.

We also note that in the case where $f$ has a zero at $z=z_0$
the two linearly independent solutions are actually 
locally holomorphic.

Finally we get the following nice result for poles of second order:

\mycorollary{} \label{n2inM} {\em
If $f$ has a pole of second order at $z=z_0$ with vanishing residue,
then \bref{THREEA} has a locally meromorphic solution at $z=z_0$.}

\separate\underline{Proof:} Condition~\bref{recrel2} in this case
reduces to $q_1=-\res_{z_0}f=0$.
\QED

\newsection \label{proof21}
At this point we are able to provide the

\separate\underline{Proof of Theorem~\ref{th21}:}
The existence of a global solution to \bref{THREEA}, \bref{THREEB} 
implies the existence of linearly independent local meromorphic solutions of
\bref{FIVE}, which are given by the entries in the right column of the
matrix solution $g_-(z,\lambda)$. Because, if these entries are multiples of
each other, then by \bref{FOURB} and \bref{FOURD} the rows of $g_-$ are
linearly dependent, which contradicts $\det g_-(z,\lambda)=1$.

To prove the converse statement we first look at the local problem: 
We therefore look at a neighbourhood $V$ of an arbitrary point $z_0$,
where $f$ has at most one pole or one zero, which is located at $z=z_0$.

Let us first consider the case where $f$ has a zero of order $n$ at
$z_0$.
Denote by $y_1$ and $y_2$ the two solutions of \bref{FIVE} associated
with $r_1=n+1$ and $r_2=0$, respectively.
We have seen in section~\ref{recursionsec} 
that $y_1$ and $y_2$ are functions of
$\lambda^{-2}$. Since, by the remark after equation~\bref{SIXA}, 
$y_1, y_2\in {\cal A}$ as functions of
$\lambda$, we can write
\BEA
y_1 & = & a_0+\lambda^{-2}a_2+\ldots, \\
y_2 & = & b_0+\lambda^{-2}b_2+\ldots
\EEA
{}From \bref{FIVE} we obtain that $a_0$ and $b_0$ satisfy the
differential equation 
\BEQ
q^\dprime-{f^\prime\over f}q^\prime=0.
\EEQ
This equation has the solution
\BEQ
q^\prime=\alpha f,\kern1cm q=\alpha\int_{z_0}^z f(z^\prime)\diff 
z^\prime +\beta,
\EEQ
for some complex constants $\alpha,\beta$.
Since $y_1(z_0)=0$ and $y_2(z_0)=1$ for all $\lambda$ we get
$a_0(z_0)=0$, $b_0(z_0)=1$. Therefore
\BEQ
a_0(z)=\alpha\int_{z_0}^z f(z^\prime)\diff z^\prime,\kern1cm 
b_0=1+\gamma\int_{z_0}^z f(z^\prime)\diff z^\prime,
\EEQ
for some $\alpha,\gamma\in\Bcc$.
Thus we can find $\sigma,\tau\in\Bcc$, $\tau\neq0$, 
s.t.\ $u=\sigma y_1+\tau y_2$
is of the form 
\BEQ
u=1+u_2\lambda^{-2}+\ldots
\EEQ 
Then $u$ and $y_1$ are linearly independent and the matrix 
\BEQ
R=\tmatrix{\rho\lambda{y_1^\prime\over f}}{\rho y_1}{
\lambda{u^\prime\over f}}{u},\;\rho\in\Bcc
\EEQ
satisfies \bref{THREEA}. Here we have to choose $\rho=\rho(\lambda)$
in such a way, that $R$ has coefficients in the Banach algebra $\FA$.
This is possible by the remark after \bref{SIXA}, since the upper left
entry of $R$ is a solution of equation \bref{SIX}.
In addition we have that
\BEQ
\det R=\rho\lambda f\inv (y_1^\prime u-y_1u^\prime).
\EEQ
is independent of $z$ and an element of $\FA$. Since, up to a
factor, it is equal to the Wronskian of \bref{FIVE} w.r.t.\ the
linearly independent solutions $y_1$ and $u$, it is nonzero for all $\lambda$.
It follows, that $R$ is invertible, and the function $\det R$ is
invertible in $\FA$.
We can therefore choose $\rho$, such that $\det R=1$ for all
$\lambda$. This implies, that the upper left entry is of the form
$1+\lambda^{-2} c_2+\ldots$ as a function of $\lambda$.

We have shown, that at points in $D$, where $f$ has a zero, there
exists a local solution to \bref{THREEA} in 
$\Lambda^-_\ast\LieSL(2,\Bcc)_\sigma$.

\separate
In the case, where $f$ has a pole of order $n$, $g$ has a zero of
order $n+m$ at $z_0$, where $m\geq0$ is the zero order of $E$ at
$z_0$. Therefore, we can argue with 
\bref{SIXA} instead of \bref{FIVE} as in the previous case.
As a consequence one obtains also in this case a solution to
\bref{THREEA} in $\Lambda^-_\ast\LieSL(2,\Bcc)_\sigma$.

\separate
We now have shown, that around every point in $D$ there exists a 
local solution of \bref{THREEA} which lies in the group
$\Lambda^-_\ast\LieSL(2,\Bcc)_\sigma$. The latter is isomorphic to the
nontwisted based negative loop group $\Lambda^-_\ast\LieSL(2,\Bcc)$ 
by the restriction of the following homomorphism from 
$\Lambda\LieSL(2,\Bcc)$ to $\Lambda^-_\ast\LieSL(2,\Bcc)$: ($n\in\Bii$)
\BEA
\lambda^n \sigma_+ & \longrightarrow & \lambda^{2n+1}\sigma_+, \\
\lambda^n \sigma_- & \longrightarrow & \lambda^{2n-1}\sigma_-, \\
\lambda^n \sigma_3 & \longrightarrow & \lambda^{2n}\sigma_3.
\EEA
The matrices
\BEQ
\sigma_+=\tmatrix0100,\;\sigma_-=\tmatrix0010,\;\sigma_3=\tmatrix100{-1}
\EEQ
are Chevalley generators of $\LieSL(2,\Bcc)$.

We choose an open cover $(U_i)_{i\in I}$ of $D$ and solutions $g_i$ 
of $\bref{THREEA}$ in the sets $U_i$, which take values in
$\Lambda^-_\ast\LieSL(2,\Bcc)_\sigma$. Let
$\tg_i:D\rightarrow\Lambda_\ast^-\LieSL(2,\Bcc)$  denote the image of $g_i$
under the automorphism given above. Since the set of poles and zeroes of
the meromorphic function $f(z)$ has no cluster points, we can assume 
that there is at most one such point in each $U_i$. Furthermore
we may assume, that each $\tg_i$ is holomorphic in $U_i$, except
possibly for one pole which is a pole or zero of $f(z)$.
By the definition of the negative loop group and the remark after
equation \bref{SIXA}, 
we know that each $\tg_i$ as a function of $\lambda$ is analytic in 
$D_\infty=\{\lambda\in\Bcc|\,|\lambda|>{1\over2}\}\cup\infty$.

Next we consider functions $h_{ij}=\tg_i\tg_j\inv$, which are defined on
the intersections $U_i\cap U_j$. These functions are
independent of $z$, since
$$
\diff h_{ij}=(\dtg_i)\tg_j\inv-\tg_i\tg_j\inv(\dtg_j) \tg_j\inv=\tg_i\xi
\tg_j\inv-\tg_i\tg_j\inv\tg_j\xi\tg_j\inv=0,
$$
and they satisfy the cocycle condition $h_{ij}h_{jk}=h_{ik}$.

We set $V_i=D_\infty$, $i\in I$, then $(V_i)_{i\in I}$ is an open 
cover of $D_\infty$. Moreover, the $h_{ij}$ form a cocycle relativ to
the $V_i$.
They define a rank 2 vector bundle over the noncompact Riemann surface
$D_\infty$.
But any vector bundle over $D_\infty$ is trivial. Therefore there exist
functions $h_i(\lambda)$ on $V_i$ taking values in 
$\LieSL(2,\Bcc)$, such that
$h_ih_j\inv=\tg_i\tg_j\inv$. We may also assume that
$h_i(\lambda)\rightarrow I$ for $\lambda\rightarrow\infty$.
The relation $h_i\inv \tg_i=h_j\inv\tg_j$ defines a globally meromorphic 
function $\tg(z,\lambda)$ on $z\in D$, which is holomorphic for
$\lambda\in D_\infty$ at each $z$,
where it is finite. This shows that $\tg(z,\lambda)$ takes values in the
based negative loop group $\Lambda^-_\ast\LieSL(2,\Bcc)$, which is
isomorphic to the twisted group $\Lambda^-_\ast\LieSL(2,\Bcc)_\sigma$.
Since $\xi$ has no pole at $z=0$, this global solution is defined
and invertible at $z=0$. It can therefore be chosen to satisfy \bref{THREEB}.
\QED

\mycorollary{} {\em If \bref{TEN} is satisfied for all points, where $f$ has
a pole or a zero, then there exist two globally 
linearly independent solutions of equation \bref{FIVE}.}

\separate\underline{Proof:} If condition \bref{TEN} is satisfied at
every pole or zero of $f$, then
locally around every point in $D$, there exist two linearly
independent solutions to \bref{FIVE}. Then Theorem~\bref{th21} implies
the existence of two globally linearly independent meromorphic solutions.
\QED

\newsection
In this section we collect the results obtained in the last sections.

\mytheorem{ 1}{\em 
\begin{itemize}
\item[a)]
\bref{FIVE} has always a meromorphic solution.
\item[b)] Let $y_1$ denote a meromorphic solution of \bref{FIVE}. Then
$y_2=y_1\cdot v$ is a solution of \bref{FIVE} if and only if
$$
v^\prime=C{f\over y_1^2},\;C=\const.
$$
\item[c)] The equation \bref{FIVE} has two linearly independent
solutions in the neighbourhood of a point $z=z_0$, if and only if 
$$
\oint_{z_0}{f\over y_1^2}\diff z=\res_{z=z_0}{f\over y_1^2}=0.
$$
\item[d)] If \bref{TEN} is satisfied for $z_0$, then the linearly 
independent solutions $y_1$ and
$y_2$ can be chosen in the following way:

If $f$ has a zero of order $n$ at $z=z_0$, then both, $y_1$ and $y_2$ are 
locally holomorphic,
$y_1$ has no zero at $z_0$ and $y_2$ has a zero of order $n+1$ at $z=z_0$.

If $f$ has a pole of order $n$, then one solution is again
locally holomorphic without a zero at $z=z_0$ and the other solution
has a pole of order $n-1$ there.
\end{itemize}}

\separate
As pointed out in section~\ref{sec24}, equation \bref{TEN} is
invariant under coordinate changes, therefore Theorem~\ref{th21} implies 

\mytheorem{ 2} {\em Equation \bref{THREEA} has a (globally)
meromorphic solution $g_-(z,\lambda)$
in the twisted loop group $\Lambda^-_\ast\LieSL(2,\Bcc)_\sigma$, 
if and only if equation \bref{TEN} is 
satisfied at every pole or zero of $f(z)$.}

\separate
There is one case, in which condition \bref{TEN} is
automatically satisfied at a pole or zero $z_0$ of $f$.

\mytheorem{ 3} \label{symth} {\em If $f$ and $E$ are symmetric in $z-z_0$,
then \bref{THREEA} has a local meromorphic solution around $z=z_0$ 
regardless whether $f$ has a pole or a zero there.}

\separate\underline{Proof:} We set $w=z-z_0$.
If $f$ and $E$ are symmetric in
$w$, then with $y_1(w)$ also $\ty(w)=y_1(-w)$ is a solution of
\bref{FIVE}. Let $y_1$ be the meromorphic solution to the higher root
$r_1$ of the indicial equation (which always exists and 
is actually holomorphic around $z=z_0$). 
We know that $\ty(z)=\alpha y_1(w)+\beta y_2(w)$, where
$y_i(w)$ is a solution to \bref{FIVE} which belongs to the root $r_i$
of the indicial equation, and $\alpha$, $\beta$ are independent of
$w$. 

We will use the discussion in sections~\ref{sec23} and~\ref{sec24}.
First let $f$ have a pole at $z=z_0$, then $y_1$ is of the 
form~\bref{y1caseone}. Therefore, 
$$
\ty(w)=1+\sum_{i=1}^\infty (-1)^i a_i w^i 
$$
has neither a pole nor a zero at $z=z_0$. This together with
the fact that $\ty(0)=y_1(0)=1$ implies, that $\ty(w)=y_1(w)$.

If $f$ has a zero of even order $n$ at $z=z_0$, then $y_1(w)$ is of the form
\bref{y1casetwo}.
Therefore, 
$$
\ty(w)=-w^{n+1}+\sum_{i=1}^\infty (-1)^i a_i w^{i+n+1} 
$$
has a zero of the same order as $y_1$ at $z=z_0$. This together with
the fact that the coefficient of $w^{n+1}$ in $\ty$ is $-1$ implies
that $\ty(w)=-y_1(w)$

Thus we have $\ty^2(w)=y_1^2(w)$ in both cases and, with $y_1^2$ and $f$
being symmetric, ${f\over y_1^2}$ cannot have a residue at $z=z_0$. 
\QED

\separate
Finally, this theorem is accompanied by a result similar to 
Corollary~\ref{n2inM}:

\mycorollary{} {\em 
Let $\xi$ be as in equation \bref{MEROPOT}. In the
following two cases \bref{THREEA} has locally a meromorphic solution
around $z=z_0$:

a) $f$ has a pole at $z=z_0$ with even principal part 
and $f\cdot E$ has at most a pole of second order there.

b) $f$ has a zero at $z=z_0$ and $f\inv E$ has no pole at
$z=z_0$.}

\separate\underline{Proof:} In Case b) neither $f$ nor $f\inv E$ has a
pole at $z=z_0$. Therefore the meromorphic potential is holomorphic
around $z_0$. Around a point where the meromorphic potential is
nonsingular, \bref{THREEA} has locally a meromorphic solution.

In Case a) $f\cdot E$ having at most a pole of second order implies, 
that $E$ has at least a zero of order $n-2$. 
But then $E_0=\ldots=E_{n-3}=0$ and the condition \bref{recrel2} 
is a condition on the principal part of $f$ only. With
Theorem~\ref{symth}.3 we then get the result.
\QED

\newsection \label{normalized}
In the following sections we will further investigate the 
condition~\bref{recrel2} in the special cases $f(z)=(z-z_0)^n$ and
$f(z)=(z-z_0)^{-n}$.

\newsection {\bf Case I:} $f$ has a pole of order $n\geq2$ at $z=z_0$.\\
We want to find a solution to \bref{FIVE} of the form
\BEQ
y_2=\sum_{i\geq0}\ta_i(z-z_0)^{i-n+1},\, \ta_0\neq0.\label{ELEVEN}
\EEQ
In this section we will normalize $f=(z-z_0)^{-n}$, $n\geq2$.
This is possible without loss of generality because of the remark
in section~\ref{sec21} about the invariance of 
\bref{THREEA} and \bref{THREEB}
under coordinate changes.

In this case, with ${f^\prime\over f}=-{n\over z-z_0}$, 
the recursion relation~\bref{recrel} reads
\BEA \label{recrelhere}
\ta_1 & = & 0,\\
k(k-n+1)\ta_k & = & \lambda^{-2}\sum_{i=2}^k E_{i-2}\ta_{k-i},\,
k\geq2.
\EEA
If we choose $\ta_0$ to be independent of $\lambda$, then, 
as was already noted in section~\ref{recursionsec}, 
that the coefficients $\ta_k$ are
even polynomials in $\lambda\inv$. In particular in this case
$\deg\ta_k\leq k$ if $k$ is even, $\deg\ta_k\leq k-1$ if $k$ is
odd.

We also have the nontrivial condition~\bref{recrel2} on $E$, which in
this case becomes:
\BEQ
\sum_{i=2}^{n-1}\ta_{n-i-1}E_{i-2}=0.\label{FOURTEEN}
\EEQ
Here the $E_i$ are the Taylor coefficients of $E(z)$:
$$
E(z)=\sum_{k=0}^\infty E_k(z-z_0)^k.
$$
Clearly, $\ta_0$ and $\ta_{n-1}$ can be chosen arbitrarily. If
$\ta_0=0$, then the corresponding solution is $\ta_{n-1}y_2$, where
$y_2$ denotes a holomorphic solution of \bref{FIVE}.

For $n=2$ the condition \bref{FOURTEEN} is trivially satisfied. 
Therefore we always find a
meromorphic matrix solution of \bref{THREEA} in this case. 

\newsection
Let now $n\geq4$ and set $\ta_0=1$. 
Then \bref{FOURTEEN} has a solution of
the form 
$$
(1,0,\ta_2,\ldots,\ta_{n-2}).
$$
This equation involves all coefficients
$E_0,E_1,\ldots,E_{n-3}$ except $E_{n-4}$, which can be chosen arbitrarily.

The left hand side of \bref{FOURTEEN} is an even polynomial in
$\lambda\inv$ of degree at most $n-3$. Its constant term is $E_{n-3}$,
therefore
\BEQ
E_{n-3}=0.
\EEQ
This implies, that if $f$ has a pole of order $n\geq4$ at $z=z_0$, and $E$
has a zero of order $n-3$ there, then there will not be a meromorphic
solution to \bref{THREEA}.

In addition we have

\label{polepowers}
\mylemma{}{\em
Let $\deg\ta_k$ denote the degree of $\ta_k$ as a polynomial of
$\lambda\inv$.
Then for $0\leq k\leq n-2$
\BEQ \label{highpow}
\deg\ta_k=\left\{\begin{array}{l} k,\;\mbox{\rm for $k$ even,}\\
k-1,\;\mbox{\rm for $k$ odd.}
\end{array}\right.
\EEQ
If $k\geq2$, then the coefficient $c_k$ of the highest $\lambda\inv$ 
power in $\ta_k$ is a positive real multiple of 
$(-1)^{k-1\over2}E_1E_0^{k-3\over2}$,
if $k$ is odd, and a positive real multiple of 
$(-1)^{k\over2}E_0^{k\over2}$ if $k$ is even.}

\separate\underline{Proof:} Proof by induction. 
For $\ta_0=1$, $\ta_1=0$ equation~\bref{highpow} holds. 
If $n\geq4$, then the lemma holds for
\BEQ
\ta_2=\lambda^{-2}(6-2n)\inv E_0,
\EEQ 
and if $n\geq6$ then it holds
also for 
\BEQ
\ta_3=\lambda^{-2}(8-2n)\inv E_1.
\EEQ
We assume that we have proved the lemma for $k\leq s-1$. Then by the recursion
relation~\bref{recrelhere} we have that 
\BEQ 
\deg\ta_s \leq \deg\ta_{s-2}+2=s,
\EEQ
if $s$ is even, and
\BEQ
\deg\ta_s \leq \deg\ta_{s-2}+2=\deg\ta_{s-3}+2 = s-1,
\EEQ
if $s$ is odd. Also, if $s$ is even, the $\lambda^{-s}$ coefficient
$c_s$ of $\ta_s$ is a negative real multiple of $E_0c_{s-2}=\alpha
(-1)^{s-2\over2}E_0^{s-2\over2}$, $\alpha\in\Brr^+$, by the induction 
assumption. 
If $s$ is odd, the $\lambda^{-s+1}$ coefficient $c_s$ of $\ta_s$ is 
a negative real multiple of 
\BEQ
E_0c_{s-2}+E_1c_{s-3}=
(\alpha+\beta)(-1)^{s-3\over2}E_1E_0^{s-3\over2},\;\alpha,\beta\in\Brr^+,
\EEQ
by the induction assumption. Therefore the lemma follows also for $k=s$.
\QED

\mycorollary{} {\em
The highest nonvanishing power of $\lambda\inv$ in equation
\bref{FOURTEEN} is $c(z)\lambda^{-n+4}$, where $c(z)$ is the
$\lambda^{-n+4}$
coefficient of $\ta_{n-3}E_0+\ta_{n-4}E_1$. The coefficient 
$c(z)$ is a nonzero real multiple of $E_1E_0^{n-4\over2}$.}

\separate\underline{Proof:} 
{}From lemma~\ref{polepowers} and $n$ even it follows, that the highest 
$\lambda\inv$ power occuring on the l.h.s.\ of \bref{FOURTEEN} is
given by $\deg\ta_{n-3}=\deg\ta_{n-4}=n-4$. 
With the notation of lemma~\ref{polepowers}, the coefficient of 
$\lambda^{-n+4}$ is $c_{n-3}E_0+c_{n-4}E_1$. Since $c_{n-3}=\alpha 
E_1E_0^{n-6\over2}$ and $c_{n-4}=\beta 
E_0^{n-4\over2}$ and $\alpha,\beta$ have the same sign, 
we get that $c(z)$ is a nonvanishing
real multiple of $E_1E_0^{n-4\over2}$. \QED

\separate
This implies that, if $f$ has a pole of order $n\geq4$ at $z_0$, and 
$f(z)=(z-z_0)^{-n}$, then 
$$
E_1E_0^{n-4\over2}=0.
$$

\newsection \label{poleexamples}
We summarize the discussion above for $n=2,4,6$.

\separate $n=2$: \bref{FOURTEEN} is trivial, therefore \bref{THREEA}
always has a meromorphic solution in a neighbourhood of a second order
pole of $f$.

\noindent $n=4$: \bref{FOURTEEN} $\Leftrightarrow E_1=0$. 

\noindent $n=6$: \bref{FOURTEEN} $\Leftrightarrow E_3=E_0E_1=0$. $E_2$
can be chosen arbitrarily.

\newsection 
Let us write the recursion relation
\BEQ
-\lambda^2 k(k-n+1)\ta_k+\sum_{i=2}^kE_{i-2}\ta_{k-i}=0
\EEQ
as a matrix relation 
\BEQ
0=Aa,\, a=(\ta_0,0,\ta_2,\ldots)^\top,
\EEQ
where
\BEQ
A=\left(\begin{array}{ccccc} 0 & \\
0 & 0 \\
E_0 & 0 & \alpha_2 & \\
E_1 & E_0 & 0 & \alpha_3 & \\
\vdots & \vdots & \vdots & & \ddots
\end{array}\right)\kern-2em\hbox to 7.2pt{\Large 0}\kern2em,\;
\alpha_k=-\lambda^2k(k-n+1).
\EEQ
We know $\ta_1=0$, therefore the second column of $A$ is irrelevant.
Also $\ta_0=1$, therefore we get 
\BEQ
0=\left(\begin{array}{c} E_0 \\ E_1 \\ \vdots\end{array}\right)
+\hat{A}\hat{a},\,\hat{a}=(\ta_2,\ldots)^\top
\EEQ
where $\hA$ is the infinite matrix
\BEQ
\hat{A}=\left(\begin{array}{ccccc} \alpha_2 & \\
0 & \alpha_3 & \\
E_0 & 0 & \alpha_4 & \\
E_1 & E_0 & 0 & \alpha_5 & \\
\vdots & \vdots & \vdots & & \ddots \end{array}\right)\kern-2em 
\hbox to 7.2pt{\Large 0}\kern2em.
\EEQ
Therefore 
\BEQ
\left(\begin{array}{c} E_0 \\ E_1 \\ \vdots \\ 
E_{n-3} \\ \vdots
\end{array}\right)=-\left(\begin{array}{ccccc} \alpha_2 & \\
0 & \alpha_3 & \\
\vdots & & \ddots & \\
E_{n-5} & \ldots & 0 & \alpha_{n-1} & \\
\vdots & & & & \ddots \end{array}\right)\kern-2em\hbox to 7.2pt{\Large
0}\kern2em\hat{a}.
\EEQ 
Because of $\alpha_{n-1}=0$, $\hat{A}$ has the structure 
\BEQ
\hat{A}=\left(\begin{array}{ccccccc} \multicolumn{2}{c}{\mbox{\Large\it B}} 
& & & & \\
E_{n-5} & \ldots & E_0 & 0 & 0 & \\
\vdots & & \vdots & E_0 & 0 & \alpha_n & \\
& & & & & & \ddots \end{array}\right)\kern-2em\hbox to 7.2pt{\Large 0}\kern2em.
\EEQ
where
\BEQ
B=\left(\begin{array}{ccccc} \alpha_2 & & & & \\ 0 & \ddots & & & \\
E_0 & \ddots & \ddots & & \\ \vdots & & \ddots & \ddots & \\
E_{n-6} & \ldots & E_0 & 0 & \alpha_{n-2}
\end{array}\right)\kern-2em\hbox to 7.2pt{\Large 0}
\EEQ
is an invertible $(n-3)\times(n-3)$ matrix with 
$\det B=\prod_{k=2}^{n-2}\alpha_k\neq0$.
Moreover, $B$ involves only $E_0,\ldots, E_{n-4}$.

It follows
\BEQ
(\ta_2,\ldots, \ta_{n-2})^\top=-B\inv(E_0,\ldots, E_{n-4})^\top.
\EEQ
This is a (fairly) explicit formula for $\ta_2,\ldots,\ta_{n-2}$.
Substituting this into \bref{FOURTEEN} we obtain
\BEQ
E_{n-3}=(E_0,E_1,\ldots,E_{n-4})S_u P B\inv(E_0, E_1, \ldots,
E_{n-4})^\top=0
\EEQ
where $S_u$ is the shift $(x_1,\ldots,x_{n-3})^\top\mapsto
(x_2,\ldots,x_{n-3},0)^\top$ and $P$ is the permutation
$(x_1,\ldots,x_{n-3})^\top\mapsto(x_{n-3},x_{n-4},\ldots,x_1)^\top$.
Setting $Q=S_uPB\inv$ we split $Q$ into a symmetric and an
antisymmetric part $Q=Q_S+Q_A$.
So we get for $\underline{E}=(E_0,E_1,\ldots,E_{n-4})^\top$ the equation
\BEQ
\langle \underline{E},Q\underline{E}\rangle=
\langle\underline{E},Q_S \underline{E}\rangle=0.
\EEQ

\newsection {\bf Case II:} $f$ has a zero of order $n\geq2$ at $z=z_0$.\\
We want to find a solution to \bref{FIVE} of the form
\BEQ
y_2=\sum_{i\geq0}\ta_i(z-z_0)^i,\, \ta_0\neq0.\label{TWELVE}
\EEQ
Again, w.l.o.g.\ we normalize in the following $f=(z-z_0)^n$, $n>0$.
Equation~\bref{recrel} in this case reads
\BEA
\ta_1 & = & 0,\\
k(k-n-1)\ta_k & = & \lambda^{-2}\sum_{i=2}^k E_{i-2}\ta_{k-i},\,
k\geq2.
\EEA
In addition, the coefficients of the holomorphic function $E$ must
satisfy condition~\bref{recrel2} which takes the form
\BEQ
\sum_{i=2}^{n+1}\ta_{n-i+1}E_{i-2}=0.\label{THIRTYTHREE}
\EEQ
$\ta_0$ and $\ta_{n+1}$ can be chosen arbitrarily.
As in Case I, for $\ta_0=0$ we get the
always existing meromorphic 
solution $y_1$ with a zero of order $n+1$ at $z_0$. 

We set $\ta_0=1$. Then \bref{THIRTYTHREE} has a solution of
the form $(1,0,\ta_2,\ldots,\ta_n)$.
The coefficients $\ta_k$ are again even polynomials in $\lambda\inv$. We
have $\deg\alpha_k\leq k$, if $k$ is even, and $\deg\alpha_k\leq k-1$, if
$k$ is odd.

Thus equation \bref{THIRTYTHREE} gives a condition on the first $n-1$ 
terms of the Taylor expansion
of $E(z)$, with the exception of $E_{n-2}$, which can be chosen arbitrarily.

The left hand side of \bref{THIRTYTHREE} is an even polynomial in
$\lambda\inv$ of degree at most $n-1$. Its constant term is $E_{n-1}$,
therefore
\BEQ
E_{n-1}=0,
\EEQ
This implies, that if $g=E\cdot f\inv$ has a pole of first order at $z=z_0$,
then there is no meromorphic solution to \bref{THREEA}.

By the same arguments as in the proof of lemma \ref{polepowers}, we get

\mylemma{}{\em
Let $\deg\ta_k$ denote the degree of $\ta_k$ as a polynomial of
$\lambda\inv$.
Then for $0\leq k\leq n$
\BEQ
\deg\ta_k=\left\{\begin{array}{l} k,\;\mbox{\rm for $k$ even,}\\
k-1,\;\mbox{\rm for $k$ odd.}
\end{array}\right.
\EEQ
If $k\geq2$, then the coefficient $c_k$ of the highest $\lambda\inv$ 
power in $\ta_k$ is a positive real multiple of
$(-1)^{k-1\over2}E_1E_0^{k-3\over2}$, if $k$ is 
odd, and a positive real multiple of $(-1)^{k\over2}E_0^{k\over2}$, 
if $k$ is even.}

\separate
Again it follows like in the last section

\mycorollary{} {\em
The highest nonvanishing power of $\lambda\inv$ in equation
\bref{THIRTYTHREE} is $c(z)\lambda^{-n+4}$, where $c(z)$ is the
$\lambda^{-n+2}$
coefficient of $\ta_{n-1}E_0+\ta_{n-2}E_1$. The coefficient 
$c(z)$ is a nonvanishing real multiple of $E_1E_0^{n-2\over2}$.}

\separate
This implies that, if $f$ has a zero of order $n\geq2$ at $z_0$,
$f(z)=(z-z_0)^n$, then
$$
E_1E_0^{n-2\over2}=0.
$$

\newsection \label{lastnorm}
In the cases $n=2$ and $n=4$ these results can be summarized
as follows

\noindent $n=2$: \bref{THIRTYTHREE} $\Leftrightarrow E_1=0$. 

\noindent $n=4$: \bref{THIRTYTHREE} $\Leftrightarrow E_3=E_0E_1=0$. $E_2$
can be chosen arbitrarily.

\separate We write the recursion relation
\BEQ
-\lambda^2 k(k-n-1)\ta_k+\sum_{i=2}^kE_{i-2}\ta_{k-i}=0
\EEQ
as a matrix relation 
\BEQ
0=Aa,\, a=(\ta_0,0,\ta_2,\ldots)^\top,
\EEQ
where
\BEQ
A=\left(\begin{array}{ccccc} 0 & \\
0 & 0 & \\
E_0 & 0 & \alpha_2 \\
E_1 & E_0 & 0 & \alpha_3 \\
\vdots & \vdots & \vdots & & \ddots
\end{array}\right)\kern-2em\hbox to 7.2pt{\Large 0}\kern2em,\;
\alpha_k=-\lambda^2k(k-n-1).
\EEQ
We have $\ta_0=1$, $\ta_1=0$, $\alpha_{n+1}=0$, and by proceeding as in
Case I we get
\BEQ
(\ta_2,\ldots, \ta_n)^\top=-B\inv(E_0,\ldots, E_{n-2})^\top.
\EEQ
where $B$ is the $(n-1)\times (n-1)$ matrix  
\BEQ
B=\left(\begin{array}{ccccc} \alpha_2 & & & & \\
0 & \alpha_3 & & & \\
E_0 & 0 & \ddots & & \\
\vdots & & & \ddots & \\
E_{n-4} & E_{n-5} & \ldots & 0 & \alpha_n
\end{array}\right)\kern-2em\hbox to 7.2pt{\Large 0}\kern2em.
\EEQ
We note that $B$ is invertible with $\det
B=\prod_{j=2}^n\alpha_j\neq0$.
Moreover $B$ involves only $E_0,\ldots,E_{n-2}$.

Condition \bref{THIRTYTHREE} reads
\BEQ
E_{n-1}=(E_0,E_1,\ldots,E_{n-2})S_u P B\inv(E_0, E_1, \ldots,
E_{n-2})^\top=0
\EEQ
where $S_u$ is the shift $(x_1,\ldots,x_{n-1})^\top\mapsto
(x_2,\ldots,x_{n-1},0)^\top$ and $P$ is the permutation
$(x_1,\ldots,x_{n-1})^\top\mapsto(x_{n-1},x_{n-2},\ldots,x_1)^\top$.
Setting $Q=S_uPB\inv$ we may split $Q$ into a symmetric and an
antisymmetric part $Q=Q_S+Q_A$.
So we get for $\underline{E}=(E_0,E_1,\ldots,E_{n-2})^\top$ the equation
\BEQ
\langle \underline{E},Q\underline{E}\rangle=\langle 
\underline{E},Q_S \underline{E}\rangle=0.
\EEQ

\section{Smoothness of the Iwasawa decomposition} \label{secondpart}

\newsection \label{firstsec}
In the following we will assume that we have found a meromorphic
solution of \bref{ONE}
\BEQ
\diff
g_-=g_-\xi,\;g_-(z,\lambda)\in\Lambda^-_\ast\LieSL(2,\Bcc)_\sigma,\; 
g_-(0,\lambda)=I.
\EEQ
We are now looking for conditions on $\xi$, s.t.\ 
\begin{description}
\item[a)] the orthogonal part
$F$ of an Iwasawa decomposition of $g_-$
\BEQ \label{Iwasawa}
g_-(z,\lambda)=F(z,\lambda) g_+\inv(z,\lambda)
\EEQ
is a smooth function in $z\in D$. Here $g_+\in\Lambda^+\LieSL(2,\Bcc)_\sigma$,
$F\in\Lambda\LieSU(2)_\sigma$,
\item[b)] the map $\Phi:D\rightarrow \Brr^3$ defined by the Sym-Bobenko
formula \bref{SymBob} is an immersion. It then automatically describes
a CMC surface in $\Brr^3$ without branchpoints.
\end{description}
We would like to point out, that whenever $\xi$ can be integrated to a
meromorphic $g_-$, it is possible to split $g_-(z,\lambda)$ like in 
\bref{Iwasawa},
where $F(z,\lambda)\in\Lambda\LieSU(2)_\sigma$, 
$g_+\in\Lambda^+\LieSL(2,\Bcc)_\sigma$.
If $g_-$ has a pole at $z=z_0$, then the Sym-Bobenko formula again 
defines a smooth CMC map on a deleted neighbourhood of $z_0$, but in 
general this map will have a singularity at $z_0$.

Our main tool to classify those meromorphic potentials which belong to
CMC immersions, will be the dressing method.

We recall, that by normalization the coefficients $f$ and $g=f\inv E$
of $\xi$ are holomorphic in a neighbourhood of $z=0$.

\newsection \label{dressbasics}
Let $\FM$ denote the set of meromorphically integrable meromorphic
potentials
\BEQ
\xi=\lambda\inv\tmatrix0f{E\over f}0\dz
\EEQ
on $D$.
We also require, that for $\xi\in \FM$ the coefficient $f$ of $\xi$ does
not vanish identically, and that $f$ has only poles and zeroes of even
order.
On $\FM$ we want to define an action of the group
$\Lambda^+\LieSL(2,\Bcc)_\sigma$. 

Let $\xi\in \FM$ and denote by $g_-$ the unique solution of \bref{ONE}.
For $h_+\in\Lambda^+\LieSL(2,\Bcc)_\sigma$ we consider $h_+g_-$. In
general,
$h_+g_-$ will not be in the big cell
$\Lambda^-_\ast\LieSL(2,\Bcc)_\sigma\cdot\Lambda^+\LieSL(2,\Bcc)_\sigma$.
However, the argument in \cite[\S2]{DPW} can be applied (almost)
verbatim to our situation and we obtain

\myprop{} {\em Let $\xi\in \FM$, $h_+\in\Lambda^+\LieSL(2,\Bcc)_\sigma$ and
denote by $g_-$ the solution of \bref{ONE}. Then there exists a
discrete 
subset $S$ of $D$ and
$\hg_-(z,\lambda)\in\Lambda^-_\ast\LieSL(2,\Bcc)_\sigma$,
$\hg_+(z,\lambda)\in\Lambda^+\LieSL(2,\Bcc)_\sigma$, s.t.\
\BEQ \label{dresssplit}
h_+(\lambda)g_-(z,\lambda)=\hg_-(z,\lambda)\hg_+(z,\lambda),
\EEQ
for all $z\in D\setminus S$. Moreover, $\hg_-$ and $\hg_+$ are
meromorphic in $D$.}

\separate By differentiation we obtain $\hxi=\hg_-\inv\dhg_-$. From
the proposition above we know, that $\hxi$ is meromorphic in $D$ and
is of the form
\BEQ
\hxi=\lambda\inv\tmatrix0{\hf(z)}{\hg(z)}0\dz.
\EEQ
{}From the definition of $\hxi$ it is clear that $\hxi\in \FM$ holds. We
set
\BEQ
\hxi=h_+.\xi.
\EEQ
{}From \bref{dresssplit} and $\hxi=\hg_-\inv\dhg_-$ it follows
\BEQ
\hg_-\inv\dhg_-=-\dhg_+\cdot \hg_+\inv+\hg_+g_-\inv
\dg_-\cdot\hg_+\inv
\EEQ
or
\BEQ \label{hpx}
h_+.\xi=\hxi=-\dhg_+\cdot\hg_+\inv+\hg_+\xi\hg_+\inv.
\EEQ
It follows that $\hg=\hf\inv E$. As a consequence, the (possibly
singular) CMC surfaces associated with $\xi$ and $\hxi$ have the same
Hopf differential.

{}From the uniqueness of the Riemann-Hilbert splitting it is easy to
derive, that \bref{hpx} defines a group action of
$\Lambda^+\LieSL(2,\Bcc)_\sigma$ on $\FM$.

Finally we note

\label{dresssmooth}
\mytheorem{} {\em Let $\xi\in \FM$ and
$h_+\in\Lambda^+\LieSL(2,\Bcc)_\sigma$. Then the CMC map
associated with $\xi$ is smooth near $z=z_0$ iff the CMC map
associated with $h_+.\xi$ is smooth near $z=z_0$.}

\separate\underline{Proof:} Assume $\xi$ leads to a smooth CMC map. 
Then $g_-=Fg_+\inv$ for some smooth
$F\in\Lambda\LieSU(2)_\sigma$. Therefore
$h_+g_-=h_+Fg_+\inv=\hg_-\hg_+=\hF\hp_+$, where
$\hF\in\Lambda\LieSU(2)_\sigma$ and
$p_+\in\Lambda^+\LieSL(2,\Bcc)_\sigma$ are defined by
the last equation. Clearly, $\hF$ defines the CMC map associated
with $\hxi=h_+.\xi$.

Moreover, $h_+F=\hF\hp_+g_+$. Comparing this to the Iwasawa splitting
$h_+F=\tF\ts_+$ of $h_+F$ we see that $\hF$ is smooth---possibly up
to a $\lambda$-independent diagonal factor $k(z,\zbar)$, $\hF=\tF k$.
But the Sym-Bobenko formula yields the same immersion for $\hF$ and
$\tF$.
This proves that the CMC map associated with $\hxi=h_+.\xi$ is
also smooth. Furthermore, $\hxi=h_+\inv.\xi$, since \bref{hpx} defines a group
action. Therefore, the converse direction of the equivalence
follows by the same arguments. \QED

\mycorollary{ 1} {\em The CMC map associated with $\xi\in \FM$
is smooth near
$z=z_0$ iff there exists some $h_+\in\Lambda^+\LieSL(2,\Bcc)_\sigma$
such that the CMC map associated with $h_+.\xi$ is smooth near
$z=z_0$.}

\separate\underline{Proof:} If the CMC map associated with $\xi$ 
is smooth then we choose
$h_+= I$. On the other hand, if an $h_+$ exists, s.t.\ $h_+.\xi$
is associated with a smooth CMC map, then by
Theorem~\ref{dresssmooth} the CMC map associated to $\xi$
is also smooth.
\QED

\mycorollary{ 2} {\em Let $\xi\in \FM$ and
$h_1,\ldots,h_m\in\Lambda^+\LieSL(2,\Bcc)_\sigma$. Assume, that
$h_1.(h_2.(\ldots h_m.\xi)\ldots)$ is holomorphic near $z=z_0$. Then the
CMC map associated with $\xi$ is smooth near $z=z_0$.}

\separate\underline{Proof:} The corollary follows from
Theorem~\ref{dresssmooth} and the fact, that the dressing \bref{hpx} 
is a group action. \QED

\newsection In this section we will discuss some `basic' dressing
transformations. It is easy to verify, that
the Lie algebra $\Lambda^+\Liesl(2,\Bcc)_\sigma$ is generated by the
following three matrices: 
\BEQ
U=\tmatrix{0}{0}{\lambda}{0},\,
D=\tmatrix100{-1},\, V=\tmatrix{0}{\lambda}{0}{0}.
\EEQ
Therefore, all dressing transformations are products (and limits)
of dressing transformations associated with these three.
We will therefore investigate what the three corresponding `basic' 
transformations do.

We consider meromorphic potentials $\xi\in \FM$ with defining functions $f$
and $E$.  What we will need and compute with are the coefficients of
$g_-$. As before we write
\BEQ \label{Gform}
g_-= \tmatrix{a}{b}{c}{d},
\EEQ
with 
\BEA \label{entryform}
a & = & 1+\sum_{i\geq1}a_i\lambda^{-i}, \nonumber\\
b & = & \sum_{i\geq1}b_i\lambda^{-i}, \nonumber\\
c & = & \sum_{i\geq1}c_i\lambda^{-i}, \nonumber\\
d & = & 1+\sum_{i\geq1}d_i\lambda^{-i}.
\EEA
The function $f(z)$ is obtained by differentiation from $b_1$ and also
{}from $c_1$
\BEQ \label{ffromentry}
f(z)=\der{b_1(z)}{z}=E(z)\left(\der{c_1(z)}{z}\right)\inv.
\EEQ
It is therefore sufficient to determine the behaviour of one of these
coefficients under a dressing transformation.

Here is our general procedure:
Let $Q$ be one of the matrizes $D$, $U$ or $V$.
Set $H=\exp(tQ)$, $t\in\Bcc$. 
Then by the definition of the dressing action we need to compute the
Birkhoff splitting $Hg_-=\tg_-\tg_+$. We will do this in detail in the
following sections.

We would like to point out here though, that in all three cases $\exp(tQ)$
is of particularly simple form:
\BEA
\exp(tU) & = & \tmatrix10{t\lambda}1, \label{eUform}\\
\exp(tV) & = & \tmatrix1{t\lambda}01, \label{eVform}\\
\exp(tD) & = & \tmatrix{e^t}00{e^{-t}}. 
\EEA

\newsection $Q=D$: We note $Hg_-=Hg_-H\inv\cdot H$ and in this case we
simply have $\tg_-=Hg_-H\inv$, $\tg_+=H$. Hence the $\lambda\inv$
coefficient of the upper right corner of $\tg_-$ is given by
\BEQ
\tb_1=e^{2t}b_1.
\EEQ
By differentiation we thus obtain 
\BEQ
\tf=e^{2t}f.
\EEQ
I.e.\ the new $b_1$ and the new $f$ are just constant multiples of the
old ones. This has been used before and also been
visualized in a beautiful movie by Charlie Gunn.

\newsection $Q=U$: 
We multiply the matrix $Hg_-$ with a matrix of the form
\BEQ
G=\tmatrix{1}{0}{q\lambda}{1},
\EEQ
where we choose $q$ such that the $\lambda$ term of $Hg_-G$ vanishes.
A straightforward computation using \bref{Gform}, \bref{entryform} and
\bref{eUform} shows
\BEQ
Hg_-G=\lambda\tmatrix{0}{0}{t+q(1+tb_1)}{0}+\tmatrix{1+qb_1}{0}{0}{1+tb_1}+
\lambda\inv\tmatrix{0}{b_1}{\ast}{0},
\EEQ
therefore $q=-{t\over1+tb_1}$. This gives for the coefficient of
$\lambda^0$ in $Hg_-G$
\BEQ
\tg_0=\tmatrix{{1\over 1+tb_1}}{0}{0}{1+tb_1}.
\EEQ
Splitting off $\tg_0$ from $Hg_-G$ we obtain $Hg_-G=p_-\tg_0$,
$p_-\in\Lambda^-_\ast\LieSL(2,\Bcc)_\sigma$.
Therefore, $\tg_+$ is the product of $\tg_0$ with $G\inv$:
\BEQ
\tg_+=\tmatrix{{1\over 1+tb_1}}{0}{t\lambda}{1+tb_1}
\EEQ
and $\tg_-=Hg_-\tg_+\inv=Hg_-G\tg_0\inv$. 
The upper right corner $\tb_1$ of the
$\lambda\inv$-coefficient of $\tg_-$ is therefore
\BEQ
\tb_1=b_1(1+tb_1)\inv.
\EEQ
A differentiation using \bref{ffromentry} yields
\BEQ \label{fUtrafo}
\tf=f(1+tb_1)^{-2}.
\EEQ

\separate{\large\bf Remark:} The formula \bref{fUtrafo} is very
interesting, since, using the dressing generated by $U$, we move
{}from a constant $f$ this way to an $\tf$ with a quadratic pole.  
But this formula also shows,
that we move away from a pole of order $n$ with any arbitrarily small
$t\neq0$ and into a zero of degree $n-2$. We will exploit this in
detail later.

\newsection $Q=V$: We multiply the matrix $Hg_-$ with a matrix
of the form
\BEQ
G=\tmatrix{1}{q\lambda}{0}{1},
\EEQ
where we choose $q$ such that the $\lambda$ term of $Hg_-G$ vanishes.
Again, a straightforward computation using \bref{Gform},
\bref{entryform} and \bref{eVform} shows
\BEQ
Hg_-G=\lambda\tmatrix{0}{t+q(1+tc_1)}{0}{0}+\tmatrix{1+tc_1}{0}{0}{1+qc_1}+
\lambda\inv\tmatrix{0}{\ast}{c_1}{0},
\EEQ
therefore $q=-{t\over1+tc_1}$. This gives for the coefficient of
$\lambda^0$ in $Hg_-G$
\BEQ
\tg_0=\tmatrix{1+tc_1}{0}{0}{{1\over 1+tc_1}}.
\EEQ
We can again split off $\tg_0$ from $Hg_-G$ and obtain
$Hg_-G=p_-\tg_0$, $p_-\in\Lambda^-_\ast\LieSL(2,\Bcc)_\sigma$.
Therefore, $\tg_+$ is again the product of $\tg_0$ with $G\inv$:
\BEQ
\tg_+=\tmatrix{1+tc_1}{t\lambda}{0}{{1\over 1+tc_1}}
\EEQ
and $\tg_-=Hg_-\tg_+\inv=Hg_-G\tg_0\inv$. 
The lower left corner $\tc_1$ of the
$\lambda\inv$-coefficient of $\tg_-$ is therefore
\BEQ
\tc_1=c_1(1+tc_1)\inv.
\EEQ
A differentiation using \bref{ffromentry} gives
\BEQ \label{fVform}
\tf=f(1+tc_1)^2.
\EEQ

\separate{\large\bf Remark:}
It is only possible to add new poles of $f$ with this
transformation if $f$ has zeroes.

\newsection Writing the results above in a compact form, we obtain
with obvious notation
\BEA \label{ftransforms}
\TD(t) f & = & e^{2t}f, \label{FDTRANSFORM}\\
\TU(t) f & = & f(1+tb_1)^{-2}, \label{FUTRANSFORM}\\
\TV(t) f & = & f(1+tc_1)^2. \label{FVTRANSFORM}
\EEA
By the discussion of the dressing transformations in
section~\ref{dressbasics} we know, that every meromorphically integrable
meromorphic potential $\xi$ is transformed into another
meromorphically integrable meromorphic potential $\hxi$. Therefore, if
$f$ is the coefficient of a $\xi\in \FM$, the
transformations \bref{FDTRANSFORM}--\bref{FVTRANSFORM} on the group
level are really coming from a group action \bref{hpx} of 
$\Lambda^+\LieSL(2,\Bcc)_\sigma$. In our context it is therefore enough to
investigate the effect of the basic dressing transformations by their
action on the coefficient $f$.

\newsection \label{dressbehave} 
In this section we discuss in some more detail the effect of $\TU$ and
$\TV$ with regard to creating/annihilating poles or zeroes (note that
$\TD$ has no effect in this regard).

Consider $z_0\in D$ such that
\begin{description}
\item[(1)] $f$ has a zero of order $n$ and $E$ has a zero of order $m$
at $z_0$.

Expanding relative to $w=z-z_0$, we know 
\BEQ
b_1=\int
f(z)\dz=\beta_0+\beta_1 w^{n+1}+\ldots,\;\beta_1\neq0. 
\EEQ
Therefore,
the denominator of $\TU(t)f$ is 
\BEQ
1+tb_1=(1+t\beta_0)+t\beta_1 w^{n+1}+\ldots.
\EEQ
This implies that $1+tb_1$ does not vanish at $z_0$ if $t$ is small.
\item[(1a)] If $t$ is sufficiently small, then $\TU(t)f$ has a zero of
order $n$ at $z_0$. If $\beta_0\neq0$ and $t=-\beta_0\inv$, then
$\TU(t)f$ has a pole of order $n+2$ at $z_0$.

Next we consider $\TV(t)f$. Therefore we consider 
\BEQ
c_1=\int {E(z)\over f(z)}\dz
\EEQ
If $m\geq n$, then ${E(z)\over f(z)}$ is holomorphic at $z=z_0$ and
has a zero at $z_0$ of order $m-n$. In this case 
\BEQ
c_1=\gamma_0+\gamma_1 w^{m-n+1}+\ldots,\;\gamma_1\neq0,
\EEQ
and 
\BEQ
1+tc_1=1+t\gamma_0+t\gamma_1 w^{m-n+1}+\ldots
\EEQ
\item[(1b)] Assume $m\geq n$. If $t$ is sufficiently small, then
$\TV(t)f$ has a zero of order $n$ at $z_0$. If $\gamma_0\neq0$ and
$t=-\gamma_0\inv$, then $\TV(t)f$ has a zero of order $2m+2-n$ at
$z_0$.

Assume now $m<n$. Then ${E(z)\over f(z)}$ has a pole at $z_0$ of order
$n-m$. In this case 
\BEQ
c_1=\gamma w^{-(n-m-1)}+\ldots,\;\gamma\neq 0.
\EEQ
As a consequence, in this case for all $t\neq0$, $c_1$ has a pole of
order $n-m-1$ at $z_0$. This implies
\item[(1c)] Assume $m<n$. Then for all $0\neq t\in\Bcc$ the function
$\TV(t)f$ has a zero of order $2m+2-n$, if $2m+2-n\geq0$, or a
pole of order $n-2m-2$, if $2m+2-n<0$.

\separate{\large\bf Remark:} We note, that the case (1b) cannot occur
if $f$ and $E$ are associated with a CMC immersion around $z_0$:
If $m\geq n$, then $f$ and $E\over f$ are both holomorphic in a
neighbourhood of $z_0$. But then also $g_-$ is holomorphic there and
the splitting $g_-=Fg_+\inv$ produces smooth $F$ and $g_+$ in a
neighbourhood of $z_0$.
As a consequence, the Sym-Bobenko formula shows (see 
section~\ref{branchtheorem}), 
that $\Phi_z$ and $\Phi_\zbar$ vanish
at $z_0$, i.e.\ $\Phi$ has a branchpoint there.

\separate
Next we consider the case
\item[(2)] $f$ has a pole of order $n$ and $E$ has a zero of order $m$
at $z_0$.

In this case 
\BEQ
b_1=\beta w^{-n+1}+\ldots,\;\beta\neq0,
\EEQ
and 
\BEQ
(1+tb_1)^{-2}=(t\beta)^{-2}w^{2n-2}+\ldots
\EEQ
has a zero of order $2n-2$ at $z_0$, independent of the choice of
$t\neq0$.
\item[(2a)] For every $0\neq t\in\Bcc$, the function $\TU(t)f$ has a
zero of order $n-2$ at $z_0$.

Finally we consider again $\TV(t)f$. Here we need to consider 
\BEQ
c_1=\int{E(z)\over
f(z)}\dz=\gamma_0+\gamma_1w^{m+n+1}+\ldots,\;\gamma_1\neq0,
\EEQ
\BEQ
1+tc_1=1+t\gamma_0+t\gamma_1w^{m+n+1}+\ldots
\EEQ
\item[(2b)] If $t$ is sufficiently small, then $\TV(t)f$ has a pole of
order $n$ at $z_0$. If $\gamma_0\neq0$, and $t=-\gamma_0\inv$, then
$\TV(t)f$ has a zero of order $2m+2+n$ at $z_0$.
\end{description}

\newsection 
In the following section we discuss what poles and zeroes one can
generate and annihilate with dressing transformations. We recall that
$\FM$ denotes the set of meromorphic potentials, which can be integrated to
meromorphic $g_-$. Also, if $\xi\in \FM$, then $f$ does not vanish
identically (by convention).

\mytheorem{} {\em Let $\xi\in \FM$ and $N$ be any positive integer. Then 
\begin{description}
\item[a)] For every open subset $U\subset D$ there exists an open
subset $V_0\subset U$, s.t.\ for every $z_0\in V_0$ there exists an
element $\txi$ in the dressing orbit of $\xi$, s.t.\ the coefficient 
$\tf$ has a pole of order $2N$ at $z_0$.
\item[b)] For every open subset $U\subset D$ there exists an open
subset $V_1\subset U$, s.t.\ for every $z_0\in V_1$ there exists an
element $\hxi$ in the dressing orbit of $\xi$, s.t.\ the coefficient
$\hf$ has a zero of order $2N$ at $z_0$.
\end{description}}

\separate\underline{Proof:} 
First we show how to create poles of $f$.\\
Let $b_1$ and $c_1$ be as above, and
choose an arbitrary open subset $U_1\subseteq U$, s.t.\ $b_1(z)\neq0$
and $E(z)\neq0$ for all $z\in U_1$. This is possible, since $b_1$ is 
meromorphic and not identically zero, and we exclude the case of the
round sphere. For every $z\in
U_1$ the dressing transformation $\TU(t_1)$, $t_1=-b_1(z)\inv$, by
(\ref{dressbehave}.1a) applied to $n=0$, produces a
pole of order $2$ at $z$. We now look at the new
$c_1^{(1)}(z)=\int{E(z)\over\tf(z)}\dz$, where $\tf$ is defined by
\bref{FUTRANSFORM}. We choose a subset $U_2\subseteq U_1$,
s.t.\ $c_1^{(1)}$ has no zeroes in $U_2$. Again, this is possible,
since $c_1^{(1)}$ is nontrivial and meromorphic. For each $z\in U_2$
the successive application of $\TU(b_1(z)\inv)$ and
$\TV(c_1^{(1)}(z)\inv)$, by (\ref{dressbehave}.2b) applied to $n=2$, 
produces a zero of order $4$ at $z$. Here we also use, that
$E(z)\neq0$ on $U_1$.

We continue with this method producing higher order poles and zeroes
at every step. If necessary, we restrict ourselves to smaller open
subsets $U_n$ of $U$. For $f$ this produces poles of order $4k-2$, $k\in\Bnn$
and zeroes of order $4k$. On the other hand, starting with the
transformation $\TV$ reverses the role of poles and zeroes. For $f$ we get
poles of order $4k$ and zeroes of order $4k-2$, $k\in\Bnn$.
\QED

\mycorollary{} {\em Let $z_0\in D$, $\epsilon>0$ and $N$ be a positive
integer.
Then there exist $z_1,z_2\in D$, $\|z_1-z_0\|<\epsilon$,
$\|z_2-z_0\|<\epsilon$ and dressing transformations $T_1$ and $T_2$,
s.t.\ $f_1=T_1 f$ has a pole of order $2N$ at $z_1$ and $f_2=T_2 f$
has a zero of order $2N$ at $z_2$. }

\newsection \label{genpoles} 
In the last section we have seen, that by using dressing
transformations we can create poles and zeroes of $f$. In the next
sections, we are mainly interested in the question whether we can
---at least locally--- remove poles and zeroes of $f$.

For a $\xi\in \FM$ we will say, that a pole (resp.\ zero) $z_0$ of $f$
can be ``dressed away locally'', 
if there exists a dressing transformation $T$
s.t.\ $Tf$ is defined and nonzero at $z_0$. Similarly we define
``dressing away locally'' of poles and zeroes for $f\inv E$.

\mytheorem{} {\em Let $\xi\in \FM$ with coefficients $f$ and $f\inv E$.
If $E(z_0)\neq 0$, then every pole or zero of $f$ at $z_0$ can be
dressed away locally.}

\separate
\underline{Proof:} Poles and zeroes of $f$ are of even order.
If $E$ has no zeroes, then, by (\ref{dressbehave}.1c), $\TV(t)$
produces a pole of order $n-2$ out of a zero of order $n\geq4$ and
$\TU(t)$, by (\ref{dressbehave}.2a), produces a zero of order $n-2$ out
of a pole of order $n\geq4$. 
This holds for all $t\neq0$. Finally a pole (zero) of
order 2 is reduced by $\TU(t)$ ($\TV(t)$) to a point, 
where $f$ is finite and different from zero. \QED

\newsection \label{genpoles2}
We would like to generalize the theorem above to the case of a general
$E$. But it seems, that even in view of $\xi\in \FM$ this is not
possible. However, we can prove

\mytheorem{} {\em Let $\xi\in \FM$ and $m$ be the zero order of $E(z)$ at
$z_0\in D$. 

If $f$ has a zero of order $n$ at $z_0$, then 
this zero of $f$  can be dressed away locally, if
\BEQ \label{dressaway}
r(2m+4)-m-2\leq n\leq r(2m+4)
\EEQ
for some integer $r\geq1$.

If $f$ has a pole of 
order $n$ at $z_0$, then this pole can be dressed away locally, if
$n=2$ or
\BEQ \label{dressaway2}
r(2m+4)-m\leq n\leq r(2m+4)+2,
\EEQ for some integer $r\geq1$}

\separate\underline{Proof:} 
We proceed as above. Poles of $f\inv E$ cannot be simple, since $g_-$
is meromorphic. 
If $E$ has a zero at $z_0$ of order $m$ and $f$
has a zero of order $n=2k\geq m+2$ there, then $\TV(t)$,
by (\ref{dressbehave}.1c), will produce a
pole of order $n-2m-2$ at $z_0$. $\TU(t)$, by (\ref{dressbehave}.2a), 
produces a zero
of order $n-2$ from a pole of order $n$. So if we start with an $f$
which has a zero of order $n$ and an $E$ with a zero of order $m$,
then there are the following cases:
\begin{enumerate}
\item $n-2m-2\leq0$: We are done. We get a zero
of order $2+2m-n\geq0$ of $f$ which gives a zero of order $n-m-2\geq0$ of
$Ef\inv$ (with a zero or pole of order $0$ we always mean a point 
where $f$ is finite and nonvanishing). 

\item $n-2m-2>0$: In this case we have a pole of order 
$n-2m-2$ after transforming with $\TV(t)$. By (\ref{dressbehave}.2a),
this pole is reduced to a zero of order
$n-2m-4$ after transforming with $\TU(t)$. If $n-2m-4=0$, we are
done. Otherwise we repeat the whole procedure with $n-2m-4$ instead of
$n$. Notice, that in this case, if $n$ satisfies
condition~\bref{dressaway} for some
integer $r$, then $r\geq2$. Therefore $n-(2m+4)$ satifies the condition
for $r-1\geq1$.
\end{enumerate}
If we start with a pole of $f$ of order $n$, we apply first $\TU(t)$.
If $n=2$, then
we get by (\ref{dressbehave}.2a) a locally holomorphic $\tf=\TU(t)f$
without zero at $z_0$. If $n\geq4$, we get 
a zero of $\tf=\TU(t)f$ of order $n-2$. If the pole order $n$
satisfies \bref{dressaway2} for some integer $r\geq1$, 
then the zero order $n-2$ satisfies
\bref{dressaway} for the same integer $r$.
This case can then be treated like the case of a zero of $f$. \QED. 

\separate{\large\bf Remark:} 
It is important that in the
process of pole/zero order reduction there is no $Ef\inv$ which has
a pole of first order, because this is impossible for the coefficient of a
meromorphically integrable potential. We get the following result:

\mycorollary{} {\em \label{smoothbranch} Let $\xi$ be of the form
\bref{MEROPOT} and let
$m$ be the zero order of $E$ at a point
$z_0\in D$. If $f$ has a zero of order $n$ at $z_0$ and
$$
n=k(m+2)-1
$$ 
for some positive odd integer $k$, then $\xi\notin\FM$.
Similarly, if $f$ has a pole of order $n$ at $z_0$ and
$$
n=k(m+2)+1
$$ 
for some positive odd integer $k$, then $\xi\notin\FM$.
}

\separate
\underline{Proof:} In the first case, if $k=2r+1$, then
$n-r(2m+4)-m=1$. Therefore the procedure for pole/zero order reduction
described above leads to an $Ef\inv$ in the dressing orbit with a
simple pole. The second part reduces to the first one after applying a dressing
transformation $\TU(t)$ to $f$. \QED

\separate
It is clear, that the condition of the corollary can only
be satisfied by an $m$ which is odd, since $n$ is always even.

\newsection
In the last section we did not exclude the case of a CMC map with
branchpoints.
If we exclude branchpoints, i.e.\ if we restrict ourselves really to
CMC immersions over the domain $D$, then we get much more
restrictive conditions on the pole and zero orders of the functions
$f(z)$ and $E(z)$.

As above let $f$ and $f\inv E$ be the entries of a meromorphic potential,
which describes a CMC map \underline{without} branchpoints.

This additional feature of $\xi$ will be expressed by the fact, that
if $f$ has a zero at a point $z_0$ then $g=f\inv E$ must have a
pole there, otherwise we get a branchpoint of the CMC map (see 
Theorem~\ref{branchtheorem}).

We will proceed similar to \ref{genpoles2}. We will apply alternatingly
$\TU(t)$ and $\TV(s)$ until we have reached a potential $\hxi$, which
is holomorphic at $z_0$, where $f$ had originally a pole or a zero.

As pointed out in \ref{firstsec}, 
the CMC map associated with $\hxi$ extends
to a smooth immersion across $z_0$. 
Therefore, by Theorem~\ref{dresssmooth}, the CMC map associated
with the original potential $\xi$ also extends to a smooth immersion
across $z_0$.

Let us start with an $f$ that has a pole of order $n>0$ and an $E$ that
has a zero of order $m$ at $z=z_0$.

If we use $\TU(t)$ on $f$, by (\ref{dressbehave}.2a) we get an $f$
with a zero of 
order $n-2$ and a $g$ with a zero of order $m+2-n$ if $n\leq m+2$ or a
pole of order $n-m-2$ if $n>m+2$, since also for the transformed
functions we have $\hg=\hf\inv E$ by the remark after equation~\bref{hpx}.
If $n-2>0$ then the new meromorphic potential
belongs to a CMC immersion only if $n>m+2$. If $n-2=0$ then we have 
found a locally holomorphic potential in the dressing orbit of the
original meromorphic potential.

Otherwise we apply $\TV(t)$ to $f$. Since $m<n-2$,
(\ref{dressbehave}.1c) applies, and yields a new $f$ with a zero of
order $2m+4-n$ at $z_0$, if $2m+4-n\geq0$, or a pole of order $n-2m-4$
at $z_0$, if $2m+4-n<0$.

If $2m+4=n$, then $f(z_0)\neq 0$. In this case $g$ has a zero of
order $m$ at $z_0$. This means, that we arrived at a holomorphic
potential. 

If $2m+4>n$, then $f$ has a zero of order $2m+4-n$ at $z_0$, therefore
$g$ must have a pole at $z_0$, i.e.\ $m-(2m+4-n)=n-m-4<0$.

If $2m+4<n$, then $f$ has a pole of order $n-(2m+4)$ and $g$ has a
zero of order $n+m-(2m+4)$, no further conditions.

In the last case we apply $\TU(t)$ and obtain a zero of $f$ at $z_0$
of order $n-2-(2m+4)$. If this turns out to be zero, then $n=(2m+4)+2$
and $g$ has a zero of order $m$ at $z_0$. This gives a holomorphic
potential.

Otherwise $n-2-(2m+4)>0$ and for $g$ we get the condition
$m-n+2+(2m+4)<0$. Thus we can apply again $\TV(t)$. From
(\ref{dressbehave}.1c) we see that we need to consider
$k=2m+2-n+2+(2m+4)=-n+2(2m+4)$. This leads again to three cases: If
$k=0$, then we have arrived at a holomorphic potential. If $k>0$, then
$f$ has a zero of order $k$ and $g$ must satisfy $m-k<0$. If $k<0$,
then $f$ has a pole of order $-k$ and we can apply $\TU(t)$ and then
$\TV(t)$ again. The new $k$ is then of the form $-n+3(2m+4)$.
Therefore, choosing the integer $r$ maximal, so that $-n+r(2m+4)<0$,
we obtain the following sequence of zero and pole orders for $f$ and
$g$, $0\leq s\leq r-1$,
$$
\begin{array}{cccccccccc}
n & \!\!\stackrel{\TU(t)}{\longrightarrow}\!\! & n-2\!\! & \ldots & 
\!\!\stackrel{\TV(t)}{\longrightarrow}\!\! & \kappa(s) & 
\!\!\stackrel{\TU(t)}{\longrightarrow}\!\! & \kappa(s)-2 & 
\!\!\stackrel{\TV(t)}{\longrightarrow}\!\! & \ldots \\
n+m & \!\!\stackrel{\TU(t)}{\longrightarrow}\!\! & n-2-m & \ldots & 
\!\!\stackrel{\TV(t)}{\longrightarrow}\!\! & \kappa(s)+m 
& \!\!\stackrel{\TU(t)}{\longrightarrow}\!\! & \kappa(s)-2-m & 
\!\!\stackrel{\TV(t)}{\longrightarrow}\!\! & \ldots
\end{array}
$$
Here, for notational convenience, we have set
\BEQ
\kappa(s)=n-s(2m+4).
\EEQ
For $s\leq r-1$, by the choice of $r$, all integers in the top row are
positive.
The integers in the bottom row are positive as well. For the numbers
$n+m-s(2m+4)$ this follows from the positivity of the number above in
the top row. For the numbers $n-2-m-s(2m+4)$ this is a consequence of
the nonbranching assumption.
We note that, by the choice of $r$, we can apply $\TU$ at any rate.

There are two cases:

\separate{\bf Case I:} Transformation with $\TU(t)$ gives an
$f$ which has neither a pole nor a zero at $z_0$, 
i.e.\ $n-2-r(2m+4)=0$. In this case $g$ has a zero
of order $m$ at $z_0$ and we get a locally holomorphic potential. 

\separate{\bf Case II:} Transformation with $\TU(t)$ gives a zero of
$f$ of order $n-2-r(2m+4)>0$. Under this assumption we get a condition
ensuring that $g$ has a pole. This condition allows us to apply
(\ref{dressbehave}.1c). We see that the pole case does not occur
because of the choice of $r$. Therefore we have two subcases to
consider. We would also like to point out, that the condition on $g$,
$n-m-2-r(2m+4)>0$, implies all the conditions on the bottom row above.
Therefore, there will be no additional conditions to be considered.

\noindent{\bf Case IIa:} $n-(r+1)(2m+4)=0$ or

\noindent{\bf Case IIb:} $n-(r+1)(2m+4)<0$ 

In case IIa we get that $g$ has a pole of order $m+2$ and $f$ a zero
of order $2m+2$. This allows us to apply (\ref{dressbehave}.1c). We
see that $\TV(t)$ yields a locally holomorphic $f$ with $f(z_0)\neq0$,
whence a locally holomorphic potential.

In case IIb smoothness requires a pole of $g$, i.e.\
$n-m-2-r(2m+4)>0$.  Because of $n-(r+1)(2m+4)<0$,
(\ref{dressbehave}.1c) shows, that $\TV(t)f$ has again a zero. 
This implies again a pole of the transformed $g$:
$n+m-(r+1)(2m+4)=n-m-4-r(2m+4)<0$.  These two conditions for the order
of $g$ can only be satisfied simultaneously if $n-m-2-r(2m+4)=1$
But this gives a simple pole of $g$ in one of the dressing steps,
which is impossible for a meromorphically integrable potential. This
contradicts the assumption, that the original $f$ and $E$ described a
CMC immersion. Therefore, this case cannot occur.

As a consequence, the only possible cases are 
\BEQ \label{twoposs}
n=r(2m+4)\kern2cm\mbox{or}\kern2cm n=r(2m+4)+2,
\EEQ
for some integer number $r\geq0$.

If we start with an $f$ that has a zero of order $n>m$ at $z=z_0$, then
the dressing transformation $\TV(t)$ gives a new $\hf$. According to
(\ref{dressbehave}.1c) we have to consider
$2m+2-n$. If $n<2m+2$, we get a zero of $\hf$ and smoothness implies
that $-n+m+2>0$. If $n=2m+2$ we get a locally holomorphic potential.
The last case $n>2m+2$ gives a pole of order $n-(2m+2)>0$ of $\hf$.

In the last case we can apply what we have shown above and see that
$n-(2m+2)=r(2m+4)$ or $n-(2m+2)=r(2m+4)+2$ for some integer $r\geq1$. 
Altogether this yields for
$n$ the two possibilities
\BEQ \label{twoposs2}
n=r(2m+4)\kern2cm\mbox{or}\kern2cm n=r(2m+4)-2
\EEQ
for some integer number $r\geq1$.\\
Therefore, only $n<2m+2$
remains. In this case also $-n+m+2>0$, which is stronger than the
first inequality. From the outset we have $n>m$. These two
inequalities $n-m>0$, $n-m-2<0$ imply $n-m=1$. But this is the pole
order of the original $g$, a contradiction, which implies that this
case cannot occur.

\newsection \label{smoothwobranch} 
Altogether we have shown

\mytheorem{} {\em Let
$$
\xi=\lambda\inv\tmatrix{0}{f(z)}{\frac{E(z)}{f(z)}}{0}\dz\in\FM.
$$
Then for $\xi$ to yield a smooth CMC immersion under
the DPW construction it is necessary and sufficient that for every
$z_0\in D$, where $f$ has a pole or a zero, the following conditions
are satisfied:
Let $m$ denote the zero order of $E$ at $z_0$. 
If $f$ has a pole of order $n$ at $z_0$, then $n=2$ or for some 
integer $r\geq1$
\BEQ \label{maincond1}
n=r(2m+4)\kern2cm\mbox{or}\kern2cm n=r(2m+4)+2.
\EEQ
If $f$ has a zero of order $n$ at $z_0$, then for some integer
$r\geq1$
\BEQ \label{maincond2}
n=r(2m+4)\kern2cm\mbox{or}\kern2cm n=r(2m+4)-2.
\EEQ}

\mycorollary{ 1} {\em If $\xi\in\FM$ and $E$ has no zeroes on $D$, 
then $\xi$ yields a CMC immersion.}

\separate\underline{Proof:} Here $m=0$, whence \bref{maincond1} and
\bref{maincond2} only say, that
$f$ has only zeroes and poles of even order. But this is already part
of $\xi\in \FM$.\QED

\separate{\large\bf Remark:}
\begin{enumerate}
\item
The conditions stated in the theorem above seem to be much more 
restrictive than the integrability conditions given in the last chapter.
\item
The conditions in the Corollary~\ref{smoothbranch} cannot be
satisfied for a pair $n$, $m$, which satisfies the conditions of
Theorem~\ref{smoothwobranch}. This can easily be proved directly.
\end{enumerate}

Finally we have the following interesting result:

\mycorollary{ 2} {\em Let $\xi$ be a meromorphic potential of the form
\bref{MEROPOT}, where $E=f\cdot g$ is a holomorphic function and $f$ is
a meromorphic function without a pole at $z=0$. If furthermore $f$ has
no zeroes  and only poles of second order with vanishing residues,
then $\xi$ is associated with a CMC immersion.}

\separate\underline{Proof:} 
By Corollary~\ref{n2inM} and the assumptions on the
coefficients of $\xi$, we get that $\xi\in\FM$. Then the corollary
follows immediately from the preceding theorem. \QED

\section{Three examples} \label{examples}
\newsection 
In this chapter we construct two examples in the dressing orbit of the
cylinder,
\BEQ
\xi=\lambda\inv\tmatrix{0}{1}{1}{0}\diff z,
\EEQ
and one example in a different dressing orbit.

Consider $\xi$ as above.
In this case, because $b_1(z)=z$ has no zeroes except at $z=0$, 
we can use the transformation $\TU(t_0)$ to generate a
pole of order $2$ at the point $z=z_0\neq0$ if we choose $t_0=-{1\over
b_1(z_0)}$. 

The resulting meromorphic potential is of the form \bref{MEROPOT}
with $E(z)=1$ and
\BEQ
f(z)=\left(\frac{z_0}{z_0-z}\right)^2,\; g(z)=f(z)\inv.
\EEQ
It belongs to a smooth surface without umbilics.

In Figure~\ref{fig1} we show part of the surface associated with
$z_0={1\over4}$. 

\begin{figure}[htbp]
\caption{\label{fig1}}
\end{figure}

\newsection
Next we would like to use $\TV(t)$ for some $t$ 
to generate an additional zero at
a point $z_1$. We need to compute
\BEQ
c_1(z)=\int_0^z{E(z)\over f(z)}\dz=\frac{z^3-3z_0z^2+3z_0^2z}{3z_0^2}.
\EEQ
This function has three zeroes, one at $z=0$ and the other two at 
$z=z_0(1-e^{\pm{2\pi\over3}})$.
Except for $z_1$ being one of these points we can, by (\ref{dressbehave}.2b), 
using $\TV(t)$ to add a zero at $z_1$ if we choose 
\BEQ \label{t0choice}
t=t_0=-c_1(z_1)\inv.
\EEQ
It turns out, that we get not just one, but rather three zeroes of 
$\hf=\TV(t_0)f$, one at $z=z_1$, the other two at
$z=z_0+(z_1-z_0)(1-e^{\pm{2\pi i\over3}})$. 
A straightforward calculation using \bref{t0choice} and \bref{fVform}
shows that
$$
\hf(z)=\left({z_0((z_0-z_1)^3-(z_0-z)^3)\over
(z_0-z)((z_0-z_1)^3-z_0^3)}\right)^2.
$$

\newsection
As a last example we consider a meromorpic potential with a pole of sixth
order at some $z=z_0\neq0$. If $f(z)=(z-z_0)^{-6}$, then we know from
Theorem~\ref{smoothwobranch} that the Hopf 
differential must be either of order $0$
or $1$ at $z=z_0$. Let us take $E(z)=z-z_0$. This $E$ automatically
satisfies the meromorphic integrability condition~\bref{TEN}
as was explained in the examples in section~\ref{poleexamples}. 
Since in addition $n=6$ and $m=1$ satisfy the condition of
Theorem~\ref{smoothwobranch}, we get that this meromorphic potential
is associated with a CMC immersion.

The surface associated with this meromorphic potential for
$z_0={1\over2}$ is partially shown in figure~\ref{fig3}.
\begin{figure}[htbp]
\caption{\label{fig3}}
\end{figure}

\begin{appendix} 
\section{Appendix} 
In this appendix we will review the Sym-Bobenko formula and explain our
conventions, which differ for example from those used in \cite{Bob}.

\newsection \label{A1} 
Let $\Phi:D\rightarrow\Brr^3$ be an immersion, 
$D\subset\Bcc$ the open unit disk or the whole complex plane.
If we take $\Phi$ as a conformal chart on the surface $\Phi(D)$, and
the metric on $\Phi(D)$ as $\diff s^2=e^u(\diff x^2+\diff y^2)$,
$u:D\rightarrow\Brr$, then the natural frame, 
\BEQ
(\Phi_x,\Phi_y,N),\kern2cm 
N={\Phi_x\times\Phi_y\over|\Phi_x\times\Phi_y|}
\EEQ
associated with $\Phi$ satisfies
\BEQ
\langle\Phi_x,\Phi_x\rangle=\langle\Phi_y,\Phi_y\rangle=e^u,\kern2cm\langle\Phi_x,
\Phi_y\rangle=0.
\EEQ
This shows that $(\Phi_x,\Phi_y,N)$ is an orthogonal frame and
\BEQ
\FU=(e^{-{u\over2}}\Phi_x,e^{-{u\over2}}\Phi_y,N)
\EEQ 
is an orthogonal matrix. By possibly rotating the surface, we can assume
$\FU(0,0)=I$. We will use this normalization from now on.

We choose complex coordinates $z=x+iy$, $\zbar=x-iy$. Then
$\partial_z={1\over2}(\partial_x-i\partial_y)$,
$\partial_\zbar={1\over2}(\partial_x+i\partial_y)$ and
\BEQ
\langle\Phi_z,\Phi_z\rangle=\langle\Phi_\zbar,\Phi_\zbar\rangle=0,\kern2cm
\langle\Phi_z,\Phi_\zbar\rangle={1\over2}e^u.
\EEQ
Using the definitions 
\BEQ
Q=\langle\Phi_{zz},N\rangle,\kern2cm H=2e^{-u}\langle\Phi_{z\zbar},N\rangle
\EEQ
we can write the Weingarten map (conventions as in \cite{Spivak}) as
\BEQ
II\cdot I\inv=\tmatrix{(Q+\Qbar)e^{-u}+H}{i(Q-\Qbar)e^{-u}}{i(Q-\Qbar)e^{-u}}{
-(Q+\Qbar)e^{-u}+H}
\EEQ
Therefore the mean curvature is given as ${1\over2}\tr(II\cdot
I\inv)=H$, justifying the nomenclature, and the Gau\ss\ curvature is
given by $K=\det(II\cdot I\inv)=H^2-4(Q\Qbar)e^{-2u}$.

\newsection \label{A2} 
In this section we discuss briefly the system of equations
\BEQ \label{firstsys}
\FU\inv\FU_z=A\kern1cm\mbox{\rm and} \kern1cm\FU\inv\FU_\zbar=B. 
\EEQ
This means
we need to compute $\partial_z(e^{-{u\over2}}\Phi_x)$,
$\partial_z(e^{-{u\over2}}\Phi_y)$ and $\partial_z N$, and the
corresponding expressions for $\partial_\zbar$. A somewhat tedious and
lengthy but standard computation yields $\partial_z\FU=\FU A$,
where
\BEQ
A=\left(\begin{array}{ccc} 0 & {i\over2}u_z &
-(Q+{1\over2}e^uH)e^{-{u\over2}} \\ -{i\over2}u_z & 0 &
-i(Q-{1\over2}e^uH)e^{-{u\over2}} \\ (Q+{1\over2}e^u H)e^{-{u\over2}}
& i(Q-{1\over2}e^uH)e^{-{u\over2}} & 0 \end{array}\right).
\EEQ
Then $\partial_\zbar\FU=\FU B$, where $B=\overline{A}$.

\mytheorem{} {\em The system of equations $\partial_z\FU=\FU A$ and
$\partial_\zbar\FU=\FU B$ has the compatibility condition
\BEQ
A_\zbar-B_z-[A,B]=0.
\EEQ
This is equivalent with 
\BEA \label{GaussCod1}
u_{z\zbar}+{1\over2}e^uH^2-2e^{-u}|Q|^2 & = & 0, \\
Q_\zbar={1\over2}e^uH_z.
\EEA}

\separate\underline{Proof:} Straightforward.

\newsection \label{A3} 
In the following we adopt the spinor representation for vectors in
$\Brr^3$, i.e.\ we identify $\Brr^3$ with $\Liesu(2)$.

The map $J:\Brr^3\rightarrow \Liesu(2)$ is defined by 
$J\vecr\mapsto-{i\over2}\vecr\vecsigma$, where $\vecsigma$ is the 
vector, whose components are the Pauli matrices
\BEQ
\sigma_1=\tmatrix{0}{1}{1}{0},\;\sigma_2=\tmatrix{0}{-i}{i}{0},\;\sigma_3=
\tmatrix{1}{0}{0}{-1}.
\EEQ
Therefore 
\BEQ
J(x,y,z)={1\over2}\tmatrix{-iz}{-ix-y}{-ix+y}{iz},
\EEQ
and we have
\BEQ
\langle\vecr_1,\vecr_2\rangle=-2\tr(J\vecr_1\cdot J\vecr_2)
\EEQ
and
\BEQ
J(\vecr_1\times\vecr_2)=[J\vecr_1,J\vecr_2].
\EEQ
This means, that $J$ is an isomorphism of $\Brr^3$ equipped with the
crossproduct and $\Liesu(2)$.
The formula above states, that the natural scalar product of $\Brr^3$
corresponds via $J$, up to a factor, to the (invariant) Killing
form of the semisimple Lie algebra $\Liesu(2)$.

Clearly, for every linear map $A$ of $\Brr^3$ the map $JAJ\inv$ is a
linear map of $\Liesu(2)$, and every linear map of $\Liesu(2)$ is of
this form.

For every invertible linear map of $\Brr^3$ we know
\BEQ
(A\vecr_1)\times(A\vecr_2)=(\det A)(A\inv)^\top(\vecr_1\times\vecr_2).
\EEQ
Therefore, the group of automorphisms of $\Brr^3$ equipped with the
cross product is the group
\BEQ
G=\{A\in\LieGL(3,\Brr);\;A=\det A(A\inv)^\top\}.
\EEQ
The defining equation implies $\det A=(\det A)^3(\det A)\inv$, whence
$\det A=1$ and $A=(A\inv)^\top$. Therefore,
\BEQ
G=\LieSO(3).
\EEQ
This proves that the group $\Aut(\Liesu(2))$ of automorphisms of
$\Liesu(2)$ is given by
\BEQ
\Aut(\Liesu(2))=J\LieSO(3)J\inv.
\EEQ
In particular, $\Aut(\Liesu(2))$ is connected. On the other hand,
there is a natural map $\LieSU(2)\rightarrow\Aut(\Liesu(2))$, where
$P\in\LieSU(2)$ acts as an inner automorphism $X\mapsto PXP\inv$.
Since also $\LieSU(2)$ is connected, the latter map is surjective.

\mytheorem{} {\em There exists a smooth map $P:D\rightarrow \LieSU(2)$,
such that the following diagram commutes

\separate\begin{center}
\setlength{\unitlength}{0.01in}%
\begin{picture}(290,93)(120,520)
\thicklines
\put(310,525){\vector( 1, 0){ 60}}
\put(155,540){\vector( 3, 2){ 75}}
\put(310,595){\vector( 2,-1){100}}
\put(150,525){\vector( 1, 0){ 80}}
\put(120,520){\makebox(0,0)[lb]{$D$}}
\put(380,520){\makebox(0,0)[lb]{$\Aut(\Liesu(2))$}}
\put(160,580){\makebox(0,0)[lb]{$P$}}
\put(190,535){\makebox(0,0)[lb]{$\FU$}}
\put(240,600){\makebox(0,0)[lb]{$\LieSU(2)$}}
\put(240,520){\makebox(0,0)[lb]{$\LieSO(3)$}}
\end{picture}
\end{center}

\separate where the unlabeled maps have been defined above. Moreover, $P$ is
unique up to a sign.}

\separate\underline{Proof:} Denote by $\hat{\FU}$ the composition of the
two maps in the bottom row. Since $D$ and $\LieSU(2)$ are simply
connected, this map factors through $\LieSU(2)$.\QED

\mycorollary{} {\em For every $X\in\Liesu(2)$ there holds 
$$
(J\FU J\inv)(X)=PXP\inv.
$$}

\newsection We extend $J$ to a $\Bcc$-linear map from $\Bcc^3$ to
$\Bcc\otimes\Liesu(2)\cong\Liesl(2,\Bcc)$. Then the complex
conjugation $\vecr\mapsto\overline{\vecr}$ of $\Bcc^3$ satisfies
\BEQ
J\overline{\vecr}=-\overline{(J\vecr)}^\top.
\EEQ
Using this and Theorem~\ref{A2} we show

\myprop{} {\em
\begin{description}
\item[(1)] $(J\FU_z J\inv)(X) = [P_zP\inv,PXP\inv],$
\item[(2)] $(J\FU_\zbar J\inv)(X) = [P_\zbar P\inv,PXP\inv],$
\item[(3)] $(J(\FU\inv\FU_z)J\inv)(X) = [P\inv P_z,X],$
\item[(4)] $(J(\FU\inv\FU_\zbar) J\inv)(X) = [P\inv P_\zbar,X].$
\end{description}}

\separate\underline{Proof:} These relations follow by a
straightforward computation from Corollary~\ref{A3}.\QED

\separate The proposition above suggests to translate from a linear
$3\times3$-system for $\FU$ to a linear $2\times2$-system for $P$.

\mytheorem{} {\em For the map $P:D\rightarrow\LieSU(2)$ we have
\begin{description}
\item[(1)] 
$$
U=P\inv P_z=\tmatrix{-{1\over4}u_z}{Qe^{-{u\over2}}}{
-{1\over2}e^{u\over2}H}{{1\over4}u_z},
$$
\item[(2)] 
$$
V=P\inv P_\zbar=\tmatrix{{1\over4}u_\zbar}{{1\over2}He^{u\over2}}{
-\Qbar e^{-{u\over2}}}{-{1\over4}u_\zbar},
$$
\item[(3)] $P(0,0)=I$.
\end{description}
Moreover, the compatibility condition for (1) and (2) is the same as
the one for \bref{firstsys}.}

\separate\underline{Proof:} In section \ref{A2} we have stated the
form of $\FU\inv\FU_z$ and $\FU\inv\FU_\zbar$. This
makes it straightforward to compute $(J(\FU\inv\FU_z)J\inv)(X)$
and $(J(\FU\inv\FU_\zbar)J\inv)(X)$. Comparing the resulting
expressions with $[R,X]$ produces (1) and (2). The equation in (3) 
follows from our normalization $\FU(0,0)=I$. The last
claim is clear, since $J$ is an isomorphism.\QED

\newsection \label{A5}
Similar to Corollary~\ref{A3} we state 
\BEQ
J(\FU e_j)=(J\FU J\inv)(Je_j).
\EEQ
We note that the definition of $J$ and $\sigma_j$ implies 
\BEQ
Je_j=-{i\over2}\sigma_j.
\EEQ
This shows
\BEA
J(e^{-{u\over2}}\Phi_x) & = & -{i\over2}P\sigma_1P\inv,\\
J(e^{-{u\over2}}\Phi_y) & = & -{i\over2}P\sigma_2P\inv,\\
J(N) & = & -{i\over2}P\sigma_3P\inv.
\EEA
As a consequence we obtain 
\BEA
J\Phi_z & = & -{i\over2}e^{u\over2}P\sigma_-P\inv, \\
J\Phi_\zbar & = & -{i\over2}e^{u\over2}P\sigma_+P\inv,
\EEA
where 
\BEQ
\sigma_+=\tmatrix0100\kern1cm\mbox{\rm and}\kern1cm\sigma_-=\tmatrix0010.
\EEQ

\newsection \label{A6}
We look at the compatibility conditions stated in
Theorem~\ref{A2} a bit more closely. We note that
if $H$ is constant, then $Q$ is holomorphic. $Q\dz^2$ is therefore a
holomorphic differential for CMC immersions, the Hopf differential.
Moreover, for differentiable surfaces condition~\bref{GaussCod1} is 
automatically satisfied.

By the standard theory of Bianchi for CMC surfaces, we know, that to
every CMC immersion $\Phi$, there belongs an associated family of CMC
surfaces.
These are classically defined by the substitution $Q\mapsto \lambda^{-2}
Q$, where $\lambda=e^{it}$ is on the unit circle in $\Bcc$.

If a function $u$ solves \bref{GaussCod1} for $Q$, then the same $u$
solves it for $\lambda^{-2} Q$. Therefore, to the same solution $u$ there
belongs a one parameter family of surfaces, which are determined by
the integral $P(\lambda)$, which depends on $\lambda$ through $Q$.

We would also like to point out that in the DPW method we work with 
matrices in a twisted loop group.
The introduction of $\lambda$ in the classical way doesn't give
elements in this group. We actually have
\BEA
U(\lambda) & = & \tmatrix{-{1\over4}u_z}{
\lambda^{-2}Q e^{-{u\over2}}}{-{1\over2}He^{u\over2}}{{1\over4}u_z},\\
V(\lambda) & = & \tmatrix{{1\over4}u_\zbar}{
{1\over2}He^{u\over2}}{-\lambda^2\Qbar e^{-{u\over2}}}{-{1\over4}u_\zbar}.
\EEA
By conjugation with the $z$-independent matrix
\BEQ \label{GDEF}
G=i\tmatrix0{\lambda^{-{1\over2}}}{\lambda^{1\over2}}0
\EEQ
we get elements $\hU(\lambda)$ and $\hV(\lambda)$ 
in the twisted loop algebra $\Lambda\Liesl(2,\Bcc)$:
$P\mapsto F=G\inv PG$, and
\BEA
U\mapsto \hU=G\inv UG & = & \tmatrix{{1\over4}u_z}{
-{1\over2}\lambda\inv He^{u\over2}}{
\lambda\inv Qe^{-{u\over2}}}{-{1\over4}u_z}, \label{UHAT} \\
V\mapsto \hV=G\inv VG & = &
\tmatrix{-{1\over4}u_\zbar}{-\lambda\Qbar e^{-{u\over2}}}{
{1\over2}\lambda He^{u\over2}}{{1\over4}u_\zbar}. \label{VHAT}
\EEA
We note that $F(0,0,\lambda)=I$, therefore $F$ lies in the twisted loop
group $\Lambda\LieSU(2)$ 

\newsection \label{Sym}
By the transformations carried out so far, the frame
$\FU$ has been translated into the ``extended lift'' $F$. In this setting
the Sym-Bobenko formula provides an easy way to derive the immersion
$\Phi_\lambda$ from the extended lift $F(\lambda)$:

\mytheorem{} {\em
Let $\Phi_\lambda:D\rightarrow \Brr^3$, $D$ as above, be an associated
family of CMC immersions with mean curvature $H$. Let $G(\lambda)$ be
defined by \bref{GDEF} and $F(\lambda)$ be given by the definitions 
in section~\ref{A6}, then
$$
G(\lambda)\inv J(\Phi_\lambda)G(\lambda)=
-{1\over 2H}\left(\frac{\partial F}{\partial t}
F\inv+{i\over2}F\sigma_3 F\inv\right)+C,\;\lambda=e^{it},
$$
where, for each $\lambda\in S^1$, $C=C(\lambda)$ is a $z$-independent 
translation of the whole surface.}

\separate\underline{Proof:} Let ``$\cdot$'' denote differentiation w.r.t.\
$t$.
Then 
\BEA
\lefteqn{\partial_z\left(-{1\over2H}\left(\dot{F}F\inv
+{i\over2}F\sigma_3F\inv\right)\right)=} & &
\nonumber\\
& = & -{1\over2H}(\dot{F}_zF\inv-\dot{F} F\inv F_z F\inv
+{i\over 2}F_z\sigma_3F\inv-{i\over 2}F \sigma_3 F\inv F_z F\inv) \nonumber\\
& = & -{1\over2H}\Ad F(\dot{\hat{U}}+{i\over2}[\hat{U},\sigma_3]) \nonumber \\
& = & -{1\over2}\Ad F\left(\tmatrix{0}{
{i\over2}\lambda\inv e^{u\over2}}{-i{Q\over H}\lambda\inv e^{-{u\over2}}}{0}+
{i\over2}\tmatrix{0}{\lambda\inv e^{u\over2}}{
{2\over H}\lambda\inv Q e^{-{u\over2}}}{0}\right) \nonumber \\
& = & -{i\over2}e^{u\over2}\lambda\inv\Ad F\sigma_+. \nonumber
\EEA
On the other hand from section~\ref{A5} we know
\BEQ \label{JPHIZ}
J((\Phi_\lambda)_z)=-{i\over2}e^{u\over2}P\sigma_-P\inv=-{i\over2}e^{u\over2}
GFG\inv\sigma_-GF\inv G\inv=
-{i\over2}e^{u\over2}\lambda\inv G(\Ad F\sigma_+)G\inv.
\EEQ
Similarly we get
\BEQ \label{JPHIZBAR}
G(\lambda)\inv J((\Phi_\lambda)_\zbar)G(\lambda)=
-{i\over2}e^{u\over2}\lambda\Ad F\sigma_-=
\partial_\zbar\left(-{1\over2H}\left(\dot{F}F\inv
+{i\over2}F\sigma_3F\inv\right)\right).
\EEQ
This proves that the two sides of the claim can only differ by a
constant, which amounts to a translation of the whole surface. \QED

\separate{\large\bf Remark:} We actually have, that if we define 
$\widetilde{\Phi}_\lambda:D\rightarrow\Brr^3$ by
\BEQ
J(\widetilde{\Phi}_\lambda)=-{1\over 2H}\left(\frac{\partial
F}{\partial t}F\inv+{i\over2}F\sigma_3F\inv\right),
\EEQ
then $\widetilde{\Phi}_\lambda(z=0)={1\over 2H}e_3$. The maps
$\widetilde{\Phi}_\lambda$ also define an associated family of CMC
immersions, which is related to the family $\Phi_\lambda$ by a
$\lambda$-dependent Euclidean transformation in $\Brr^3$.

\newsection \label{branchtheorem}
In the DPW approach $F$ is obtained by an Iwasawa decomposition of
$g_-$. In this case, by \cite[Lemma 4.2]{DPW}, we get, that the
derivative of the Gau\ss\ map $F\inv\dF$ has the form
\BEQ \label{DPWform}
F\inv\dF=\lambda\inv\alpha_1^\prime\dz+\alpha_0+\lambda\alpha_1^\dprime\dzbar,
\EEQ
where $\alpha_0$ is a one form with values in the diagonal elements
in $\Liesl(2,\Bcc)$ and
$\alpha_1^\prime\dz$, $\alpha_1^\dprime\dzbar$ are the holomorphic and
antiholomorphic part, respectively, of a one form taking values in the
off-diagonal matrices in $\Liesl(2,\Bcc)$.

Now we prove that $\Phi$ fails to be an immersion, where the
meromorphic potential $\xi$ is holomorphic and the upper right 
corner $f(z)$ of $\xi$ has a zero.

\mytheorem{} {\em If a CMC immersion is obtained by the
DPW procedure
starting from a meromorphic potential $\xi$ of the form \bref{MEROPOT}, then
the upper right entry of $\xi$ cannot have a zero at a point, where $\xi$ is
defined.}

\separate\underline{Proof:} 
Let $\xi$ be a meromorphic potential and assume, that it is associated
with a smooth CMC map. Then the canonical moving frame
$\hF\in\Lambda\LieSU(2)_\sigma$ obtained via Iwasawa decomposition
satisfies $g_-=\hF\hg_+\inv$, $\hg_-\inv\dhg_-=\xi$,
$\hg_-(0,0,\lambda)=I$ and $\hg_+\in\Lambda^+\LieSL(2,\Bcc)_\sigma$. 
On the other hand, the decompositions used in
the body of this paper are obtained as follows: one integrates
$g_-^\prime=g_-\xi$, $g_-(0,0,\lambda)=I$ and decomposes $g_-=Fg_+\inv$,
where $F\in\Lambda\LieSU(2)_\sigma$ and
$g_+\in\Lambda^+\LieSL(2,\Bcc)_\sigma$. 

It is easy to see that $\hF=Fk$ holds, where
$k=k(z,\zbar)$ is a $\lambda$-independent diagonal matrix. If $\xi$ is
defined at $z_0$, then all occuring matrices are smooth in a
neighbourhood of $z_0$. We can therefore assume $\hF=F$. This allows
us to use the (proof of the) Sym-Bobenko formula.

We will actually prove, that the derivative
$\diff\Phi(z_0)$ vanishes if $f(z_0)=0$.
{}From the proof of the Sym-Bobenko formula above we know
\BEA
J((\Phi_\lambda)_z)=-{i\over2}e^{u\over2}\lambda\inv\Ad F\sigma_+.
\EEA
Here $s=-{1\over2}\lambda\inv He^{u\over2}$ is the upper right entry 
of $\alpha_1^\prime$, which can be read off equation~\bref{UHAT}. On the other
hand, because of $F=g_-g_+$, we get that the upper right entry of
$F\inv F_z$ is of the form 
\BEQ \label{secke}
\lambda\inv w(z,\zbar)^{-2}f(z)+\lambda^1\ldots, 
\EEQ
where $w(z,\zbar)$ is the coefficient of $\lambda^0$ in
the upper diagonal matrix entry of $g_+$. If $\xi$ is defined at
$z_0$, then $g_-$, $g_+$ and $F$ are smooth at $z_0$.
Therefore, with $w(z,\zbar)\neq0$, $(\Phi_\lambda)_z$ vanishes 
at $z_0$ if and only if $f$ does. Similarly we get that
$(\Phi_\lambda)_\zbar$ vanishes at $z_0$ if and only if $f$ does,
which proves the theorem. \QED

\separate{\large\bf Remark:}
The proof above shows that
\BEQ \label{w1}
w^{-2}f=-{1\over2}He^{u\over2}.
\EEQ
The lower left entry of $F\inv F_z$ is 
\BEQ
\lambda\inv w^2{E\over f}+\lambda^1\ldots
\EEQ
which gives, by equation \bref{UHAT},
\BEQ
w^2{E\over f}=Qe^{-{u\over2}}.
\EEQ
Using this and \bref{w1} we get $E=-{1\over2}QH$.
By comparing this with the definition of the Hopf differential in 
\cite[(4.39)]{DPW}, let us call it $\tQ$, and using \cite[(4.41)]{DPW}
we see that with our definitions $\tQ=\langle N_z,N_z\rangle=-4E=2HQ$. 
For $H=-2$ we have $Q=E$. Of course,
sign and absolute value of $H$ can be adjusted by the choice of the 
normal vector and/or a scaling of the metric in $\Brr^3$.

As another consequence we obtain a result emerging in a discussion of
Burstall and Pedit with the authors.

\mycorollary{} {\em Let $\xi$ be a meromorphic potential of a CMC surface,
and $g_-$ as above. Let $g_-=\hF\hg_+\inv$ denote the Iwasawa
decomposition of $g_-$, and $g_-=Fg_+\inv$ any decomposition of
$g_-$ with $g_+\in\Lambda^+\LieSL(2,\Bcc)_\sigma$ and
$F\in\Lambda\LieSU(2)_\sigma$ smooth. Furthermore, let 
$((g_+)_0)_{11}=ae^{ib}$ 
be the polar decomposition of the $\lambda^0$-coefficient of $\hg_+$.
Then $e^{-2ib}f$ is real.}
\end{appendix}

\end{document}